\documentclass[11pt]{article}
\usepackage{latexsym, amssymb, amsmath, amscd, amsfonts ,texdraw, mathrsfs, amsthm}
\usepackage{colordvi}

\newtheorem{thm}{Theorem}[section]

\newtheorem{cor}[thm]{Corollary}
\newtheorem{prop}[thm]{Proposition}

\theoremstyle{definition}
\newtheorem{defi}[thm]{Definition}
\newtheorem{ex}{Example}

\def\ls{large Schr\"{o}der }
\def\s{Schr\"{o}der }
\def\nlp{noncrossing linked partition }
\def\gpi{{\mathcal{G}_\pi}}
\def\ghpi{\mathcal{G}_{\hat \pi}}
\def\goto{\rightarrow}
\def\hp{\hat \pi}
\def\ml{\mathcal{L}}

\begin{document}

\title{Linked Partitions and  Linked Cycles}

\author{
William Y. C. Chen$^{1,4}$, Susan Y. J. Wu$^2$ and Catherine H. Yan$^{3,5}$ 
 \vspace{.2cm} \\
$^{1,2,3}$ Center for Combinatorics, LPMC\\
Nankai University, Tianjin 300071, P. R. China
\vspace{.3cm} \\
$^3$Department of Mathematics\\
Texas A\&M University, College Station, TX 77843 \vspace{.2cm} \\
E-mail:$^1$chen@nankai.edu.cn, $^2$wuyijun@126.com,
$^3$cyan@math.tamu.edu}
\date{
 Dedicated to Professor Zhe-Xian Wan on the Occasion of
His Eightieth Birthday}

\maketitle

\begin{abstract}
The notion of noncrossing linked partition arose from the study of certain
transforms in free probability theory. It is known that the number of
noncrossing linked partitions of $[n+1]$ is equal to the $n$-th large \s
number $r_n$, which counts the number of \s paths.
In this paper  we give a bijective proof of this
result. Then we introduce the structures of linked partitions and linked
cycles.  We present various combinatorial properties of noncrossing
linked partitions, linked partitions, and linked cycles, and
connect them  to other combinatorial structures and results,
including increasing trees, partial matchings,
$k$-Stirling numbers of the second kind, and the symmetry between
crossings and nestings over certain linear graphs.
\end{abstract}



\footnotetext[4]{The first author was supported by the 973 Project on
Mathematical Mechanization, the Ministry of Education, the Ministry
of Science and Technology, and the National Science Foundation of
China.}
\footnotetext[5]{The third author was supported in part by NSF grant \#DMS-0245526.}

\section{Introduction}
One of the most important combinatorial structures is a \emph{partition of a
finite set} $N$, that is,  a collection $\pi=\{B_1, B_2,
\ldots, B_k\}$ of subsets of $N$ such that
(i) $B_i\neq \emptyset$  for each $i$;
(ii) $B_i\cap B_j=\emptyset$ if $i\neq j$, and
(iii) $B_1 \cup B_2 \cup\cdots\cup B_k = N.$
Each element $B_i$ is called a \emph{block} of $\pi$.

Let $N=[n]$, the set of integers $\{1, 2, ...,n\}$, and $B_i$, $B_j$ be two
blocks of a partition $\pi$ of $[n]$.
We say that $B_i$ and $B_j$ are {\it crossing}
if there exist $a$, $c\in B_i$ and $b$, $d\in B_j$ with
$a<b<c<d$. Otherwise, we say $B_i$ and $B_j$ are {\it
noncrossing}. A \emph{noncrossing partition} $\sigma$ is a partition
of $[n]$ whose blocks are pairwise noncrossing.

Given partitions $\pi$ and $\sigma$ of $[n]$ we say that $\pi <
\sigma$ if each block of $\pi$ is contained in a block of $\sigma$.
This ordering defines a lattice
on the set of all partitions of $[n]$,  which is
called the partition lattice.
When restricted
to the set of noncrossing partitions on $[n]$, it is called the
\emph{noncrossing partition lattice} and denoted $NC_n$.
The noncrossing partition lattice is a combinatorial structure
that occurs in a diverse list of mathematical areas, including, for
example, combinatorics, noncommutative probability, low-dimensional topology
and geometric group theory.  A nice expository article on the
subject is given in \cite{Mccammond}.

Recently, in studying the unsymmetrized $T$-transform in the content
of free probability theory,  Dykema  introduced a new
combinatorial structure,  the
\emph{noncrossing linked partition} \cite{Dykema}, which can be viewed as
a noncrossing partition with possible some links with restricted nature
drawn between certain blocks of the partition.
Dykema described two natural partial orderings on the set of
noncrossing linked partitions of $[n]$, and compared it with the
noncrossing partition lattice $NC_n$. In particular, he obtained
the generating function for the number of noncrossing linked partitions
via transforms in free probability theory. It follows
that the cardinality of noncrossing linked partitions of $[n+1]$ is
equal to the $n$-th large Schr\"oder number $r_n$, which counts the
number of \s paths of length $n$. A \s path of length $n$ is a
lattice path from $(0,0)$ to $(n,n)$
consisting of steps East $(1,0)$, North $(0,1)$ and Northeast
$(1,1)$, and never lying under the line $y=x$.
The first few terms of the \ls numbers are
$1,2,6,22,90,394,1806\ldots$.
It is the sequence  A006318 in the  database \emph{On-line Encyclopedia of
Integer Sequences} (OEIS) \cite{Sloane}.

The restricted link between blocks proposed by Dykema is as follows.
Let $E$ and $F$
be two finite  subsets of integers.
We say that $E$ and $F$ are {\it nearly disjoint} if for every
$i\in E\cap F$, one of the following holds:
\begin{itemize}
\item[{\bf a.}] $i=\min(E)$, $|E|>1$ and $i\neq\min(F)$, or
\item[{\bf b.}] $i=\min(F)$, $|F|>1$ and $i\neq\min(E)$.
\end{itemize}

\begin{defi}
A {\it linked partition} of $[n]$ is a set $\pi$ of
nonempty subsets of $[n]$ whose  union  is  $[n]$ and
any two distinct elements of $\pi$ are nearly
disjoint.  It is a {\it noncrossing linked partition} if
in addition, any two distinct elements of $\pi$ are noncrossing.
\end{defi}

Denote by  $LP(n)$  and $NCL(n)$) the set of all linked
partitions  and noncrossing linked partitions  of $[n]$ respectively.
As before,
an element of $\pi$ is called a {\it block} of $\pi$.
In the present paper we study the combinatorial properties of linked
partitions.
Section 2 is devoted to the  noncrossing linked partitions.
We construct a bijection between the set of
noncrossing linked partitions of the set $[n+1]$ and
the set of Schr\"oder paths of length $n$, and derive various
generating functions for noncrossing linked partitions.
In Section 3 we discuss the set $LP(n)$ of all linked partition of $[n]$.
We show that $LP(n)$ is in one-to-one correspondence with the set of
increasing trees on $n+1$ labeled vertices, and derive properties of
linked partitions from those of increasing trees.
We also define two statistics for a linked partition ,
 2-crossing and 2-nesting,
and show that these two statistics are equally
distributed over all linked partitions with the same lefthand and
righthand
endpoints.
Then we propose a notion of \emph{linked cycles},
which is a linked partition equipped with a cycle structure  on each of
its block. We describe two graphic representations of linked cycles,
 give the enumeration of the linked cycles on $[n]$,
and study certain statistics over
linked cycles. In particular, we show that there are two symmetric
joint generating functions  over all linked cycles on $[n]$: one for
2-crossings and 2-nestings, and the other for the crossing number and
the nesting number. This is the content of Section 4.

\section{Noncrossing linked partitions}

This section  studies the combinatorial properties of noncrossing
linked partitions. Let $f_n=|NCL(n)|$, the number of noncrossing
linked
partitions of $[n]$.
We establish a  recurrence for the sequence $f_n$, which leads to the
generating function. Then we give a bijective proof of the
identity $f_{n+1}=r_n$, where $r_n$ is the $n$-th large
Schr\"oder number.
Using the  bijection,  we derive various enumerative results for
noncrossing linked partitions.

The following basic properties of noncrossing linked partitions were observed
in \cite[Remark 5.4]{Dykema}.

\noindent {\bf Property}. Let $\pi \in NCL(n)$.
\begin{itemize}
\item[1.] Any given element $i$ of $[n]$ belongs to either exactly one or
exactly two blocks; we will say  $i$ is singly or doubly covered by $\pi$,
accordingly.

\item[2.] The elements $1$ and $n$ are singly covered by $\pi$.

\item[3.] Any two blocks $E$, $F$ of $\pi$ have at most one element in common.
Moreover,
if $|E\cap F|=1$, then both $|E|$ and $|F|$ have at least two elements.
\end{itemize}

Noncrossing linked partitions can be represented by graphs.
One such graphical representation is described in \cite{Dykema}, which
is a modification of the usual picture of a noncrossing partition.
In this representation, for $\pi \in NCL(n)$,
one lists $n$ dots in a horizontal line, and connects the $i$-th one with the
$j$-th one if and only if $i$ and $j$ are consecutive numbers in a
block of $\pi$. Here  we propose a new graphical representation,
called \emph{the linear representation}, which
plays an important role in the bijections with other combinatorial
objects.
Explicitly, for a linked partition $\pi$ of $[n]$,
list $n$ vertices  in a horizontal line with labels $1, 2, \dots, n$.
For each block $E=\{i_1,i_2,\ldots,i_k\}$ with
$i_1=\min(E)$ and $k \geq 2$,
draw an arc between $i_1$ and $i_j$ for each $j=2,\ldots, k$.
Denote an arc by $(i,j)$ if $i <j$, and call $i$ the lefthand
endpoint, $j$ the righthand endpoint. In drawing the graph
we always put the arc $(i,j)$ above $(i,k)$ if $j > k$.
Denoted by $\gpi$ this linear representation.
It is easy to check that a linked partition is noncrossing if and only
if there are no two crossing edges in $\gpi$.

\begin{ex}
Figure \ref{f1} shows the linear representations of all (noncrossing) linked partitions in
 $NCL(3)$.

\begin{figure}[ht]
\centertexdraw{
\drawdim mm

\linewd 0.2
\move(0 0) \rlvec (0 32) \rlvec(80 0) \rlvec(0 -32) \rlvec(-80 0)
\move(0 8) \rlvec (80 0)
\move(0 16) \rlvec (80 0)
\move(24 0) \rlvec (0 32)
\move(63 0) \rlvec (0 32)

\linewd 0.4
\move(0 24) \rlvec (80 0)
\move(40 0) \rlvec (0 32)

\linewd 0.2

\move(27 18) \clvec (29 23)(35 23)(37 18)
\move(27 18) \clvec (28.5 20)(30.5 20)(32 18)
\move(27 18) \fcir f:0 r:0.4
\move(32 18) \fcir f:0 r:0.4
\move(37 18) \fcir f:0 r:0.4

\move(27 10) \clvec (28.5 13)(30.5 13)(32 10)
\move(32 10) \clvec (33.5 13)(35.5 13)(37 10)
\move(27 10) \fcir f:0 r:0.4
\move(32 10) \fcir f:0 r:0.4
\move(37 10) \fcir f:0 r:0.4

\move(32 2) \clvec (33.5 5)(35.5 5)(37 2)
\move(27 2) \fcir f:0 r:0.4
\move(32 2) \fcir f:0 r:0.4
\move(37 2) \fcir f:0 r:0.4

\move(13 28) \textref h:R v:C \htext{$\pi$}
\move(19 20) \textref h:R v:C \htext{$\{1,2,3\}$}
\move(22 12) \textref h:R v:C \htext{$\{1,2\}\{2,3\}$}
\move(20 4) \textref h:R v:C \htext{$\{1\}\{2,3\}$}
\move(34 28) \textref h:R v:C \htext{$\mathcal{G}_{\pi}$}

\move(66 18) \clvec (67.5 21)(69.5 21)(71 18)
\move(66 18) \fcir f:0 r:0.4
\move(71 18) \fcir f:0 r:0.4
\move(76 18) \fcir f:0 r:0.4

\move(66 10) \clvec (68 15)(74 15)(76 10)
\move(66 10) \fcir f:0 r:0.4
\move(71 10) \fcir f:0 r:0.4
\move(76 10) \fcir f:0 r:0.4

\move(66 2) \fcir f:0 r:0.4
\move(71 2) \fcir f:0 r:0.4
\move(76 2) \fcir f:0 r:0.4

\move(53 28) \textref h:R v:C \htext{$\pi$}
\move(60 20) \textref h:R v:C \htext{$\{1,2\}\{3\}$}
\move(60 12) \textref h:R v:C \htext{$\{1,3\}\{2\}$}
\move(61 4) \textref h:R v:C \htext{$\{1\}\{2\}\{3\}$}
\move(74 28) \textref h:R v:C \htext{$\mathcal{G}_{\pi}$}
}
\caption{The elements of $NCL(3)$ and their
 linear representations.\label{f1}}
\end{figure}

\end{ex}

Given  $\pi \in NCL(n)$,
for a singly covered  element $i\in [n]$,
denote by $B[i]$ the block containing $i$.
If $i$ is the minimal element of a block $B$ of $\pi$, we say that
it is a minimal element of  $\pi$.
Our first result is  a recurrence for the sequence $f_n$.

\begin{prop} \label{recurrence}
The sequence $f_n$ satisfies the recurrence
\begin{eqnarray} \label{rec-for}
f_{n+1}=f_n+f_1f_n+f_2f_{n-1}+\cdots+f_nf_1,
\end{eqnarray}
with the initial condition $f_1=1$.
\end{prop}
\begin{proof}
Clearly $f_1=1$.
Let $\pi \in NCL(n+1)$ and $i=\min(B[n+1] )$.

If $i=n+1$, then $n+1$ is a singleton block of $\pi$. There are $f_n$
noncrossing linked partitions satisfying this conditions.

If $1\leq i\leq n$, then for any two elements $a,b \in [n]$  with
$a < i <b$,  $a$ and $b$ cannot be in the same block.
Hence $\pi$ can be viewed as a union of two noncrossing linked partitions,
one of $\{1, 2, \dots, i\}$, and the other of $\{i, i+1, \dots, n+1\}$ where
$i$ and $n+1$ belong to the same block. Conversely, given a noncrossing
linked partition $\pi_1=\{B_1, \dots, B_k\}$ of
$\{1, 2, \dots, i\}$ with $i \in B_k$, and a noncrossing linked
partition $\pi_2=\{C_1, \dots, C_r\}$  of $\{i, i+1, \dots, n+1\}$
with $i$ and $n+1 \in C_1$, we can obtain a noncrossing linked partition
$\pi$ of $[n+1]$ by letting
\begin{eqnarray*}
\pi = \left\{ \begin{array}{ll}
            \pi_1 \cup \pi_2  & \text{if $B_k\neq \{i\}$}, \\
            \pi_1 \cup \pi_2 \setminus \{B_k\} & \text{if
            $B_k=\{i\}$. } \end{array} \right.
\end{eqnarray*}
Also note that a noncrossing linked partition of $\{i, \dots, n+1\}$
with $i$ and $n+1$ in the same block can be obtained uniquely from a
noncrossing linked partition of $\{i, \dots, n\}$ by adding $n+1$ to
the block containing $i$. Hence we get
$$
f_{n+1}=f_n+f_1f_n+f_2f_{n-1}+\cdots+f_nf_1,
$$
for all $n \geq 0$.
\end{proof}

Prop.~\ref{recurrence}   leads to an equation for the generating function
$F(x)=\sum_{n\geq 0} f_{n+1}x^n$.
\begin{eqnarray*}
F(x) = \sum_{n=0}^{\infty}f_{n+1}x^n
     &=&
     1+\sum_{n=1}^{\infty}f_nx^n+\sum_{n=1}^{\infty}(f_1f_n+f_2f_{n-1}+\cdots+f_nf_1,)x^n\\
     &=& 1+x\cdot F(x)+ x\cdot F(x)^2.
\end{eqnarray*}
It follows that
\begin{equation}\label{GF}
F(x)=\frac{1-x-\sqrt{1-6x+x^2}}{2x}.
\end{equation}
Formula \eqref{GF}  was first obtained by Dykema \cite{Dykema}
using transforms in free probability theory.
It is the same  as the generating function of the \ls numbers,
where the $n$-th large \s number $r_n$ counts the number of \s paths
of length $n$, i.e.,
lattice
paths from $(0,0)$ to $(n,n)$ consisting of steps $(1,0)$, $(0,1)$,
and $(1,1)$, and never lying under the line $y=x$.
Therefore
\begin{thm}[Dykema]
For every $n\geq 0$, the number of elements in $NCL(n+1)$ is equal to the
large Schr\"{o}der number $r_n$.
\end{thm}

Here we  construct a bijection between noncrossing linked partitions
of $[n+1]$ and \s paths of length $n$.
For convenience, we use $E$, $N$ and $D$ to denote
East, North and Northeast-diagonal steps, respectively.

\noindent \underline{{\bf A map $\phi$ from
  $NCL[n+1]$ to the set of
\s paths of length $n$.}} \\
Given a \nlp $\pi$ of $[n+1]$ $(n \geq 0)$, define a lattice path
 from the origin $(0,0)$ by the following steps.
\begin{description}
\item[Step 1.] Initially set $x=0$. Move $k-1$ $N$-steps if the block $B[1]$
contains $k$ elements.

In general, for $x=i>0$, if $i+1$
is the minimal element of a block $B$ of $\pi$ and $|B|=k$, move
$(k-1)$ $N$-steps.

\item[Step 2.] Move one $D$-step if $i+2$ is a singly covered minimal
element of $\pi$. Otherwise move one $E$-step.
Increase the value of $x$ by one.
(Note that the path reaches the line $x=i+1$ now).
\end{description}
Iterate Steps 1 and 2 until $x=n$.
When the process terminates, the resulting lattice path is $\phi(\pi)$.

\begin{thm}
The above defined map $\phi$ is a bijection from the set of noncrossing linked
partitions of $[n+1]$ to the set of Schr\"oder paths of length $n$.
\end{thm}
\begin{proof}
In $\pi$ each integer $i \in [n+1]$ is of one of the following types:
\begin{enumerate}
\item $i$ is a singly covered minimal element;
\item  $i$ is singly covered, but $i \neq \min(B[i])$;
\item $i$ is doubly covered. In this case assume $i$ belong to blocks
  $E $ and $F$ with $i = \min(F)$ and $j=\min(E)<i$;
\end{enumerate}
Each element $i$  of the first type except $i=1$
 contributes one $D$-step between
the lines $x=i-2$ and $x=i-1$.
Each element $i$ of the second type
contributes one $N$-step at the line $x=\min(B[i])$, and one
$E$-step between the lines $x=i-2$ and $x=i-1$.
Each element $i$ of the third type contributes
one $N$-step at the line $x=j$, and one $E$-step
between the lines $x=i-2$ and $x=i-1$.
Hence the path ends at
$(n,n)$, and at any middle stage, the number of $N$-steps is no less than that
of $E$-steps. This proves that the path $\phi(\pi)$ is a \s path of length
$n$.

To show that $\phi$ is a bijection, it is sufficient to give the inverse
map of $\phi$. Given a \s path of length $n$, for each $i=2, 3, \dots, n$,
check the segment between the lines $x=i-2$ and $x=i-1$.
If it is a $D$-step, then $i$ is a singly covered minimal element.
If it is a $E$-step, draw a line  with slope 1 which starts at the
middle point of this $E$-step, and lies between the {\s} path and the
line $x=y$.  Assume the line meets the given \s path for the first time
at $x=j < i-1$. Then $i$ belongs to a block whose minimal element is $j+1$.
We call this diagonal segment between $x=j$ and $x=i-\frac 32$ a
\emph{tunnel}  of the \s path. See Figure \ref{tunnel} for an illustration.

The resulting collection of subsets of $[n+1]$ must be pairwise noncrossing.
This is because for any tunnel whose endpoints are $A$ and $B$, where
$A$ is an $N$-step and $B$ is an $E$-step, there are an equal number
of $N$-step and $E$-steps between $A$ and $B$. Therefore for any
$E$-step between $A$ and $B$, the tunnel starting from it must end at an
$N$-step between $A$ and $B$ as well.
Also note that any element $i \in [n+1]\backslash\{1\}$ can belong to
at most two such subsets. If it happens, then both subsets have cardinality
at least two, and $i$ is the minimal element of exactly one of them. Hence
the collection of subsets obtained  forms a noncrossing linked partition.
We leave to the reader to check that this  gives  the inverse of $\phi(\pi)$.
\end{proof}

\begin{ex}
Tunnels in a \s path of length 5. \\
 The tunnel between $x=1$ and $x=3/2$ implies that $3$ is in a
  block $B$ with $\min(B)=2$.  \\
 The tunnel between $x=2$ and $x=3/2$ implies that $4$ is in a
  block  $B$ with $\min(B)=3$. \\
 The tunnel between $x=1$ and $x= 7/2$ implies that $5$ is in a
  block   $B$ with $\min(B)=2$.  \\
 The tunnel between $x=0$ and $x=9/2$ implies that $6$ is in a
  block  $B$ with $\min(B)=1$.

The corresponding noncrossing linked partition is $\pi=\{\{1,6\},
\{2,3,5\}, \{3,4\}\}$.

\begin{figure}[h,t]
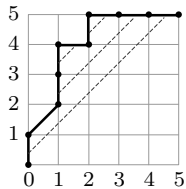

\centertexdraw{
\drawdim mm

\linewd 0.05 \setgray 0.6
\move(10 10) \rlvec(0 20)\rlvec(20 0)\rlvec(0 -20)\rlvec(-20 0)
\move(14 10) \rlvec(0 20)
\move(18 10) \rlvec(0 20)
\move(22 10) \rlvec(0 20)
\move(26 10) \rlvec(0 20)
\move(10 14) \rlvec(20 0)
\move(10 18) \rlvec(20 0)
\move(10 22) \rlvec(20 0)
\move(10 26) \rlvec(20 0)

\linewd 0.3 \setgray 0
\move(10 10) \rlvec(0 4) \rlvec(4 4) \rlvec(0 8) \rlvec(4 0) \rlvec(0 4) \rlvec(12 0)

\move(10 10) \fcir f:0 r:0.4
\move(10 14) \fcir f:0 r:0.4
\move(14 18) \fcir f:0 r:0.4
\move(14 22) \fcir f:0 r:0.4
\move(14 26) \fcir f:0 r:0.4
\move(18 26) \fcir f:0 r:0.4
\move(18 30) \fcir f:0 r:0.4
\move(22 30) \fcir f:0 r:0.4
\move(26 30) \fcir f:0 r:0.4
\move(30 30) \fcir f:0 r:0.4

\linewd 0.1 \setgray 0.2
\lpatt(0.5 0.3)
\move(10 11.5) \rlvec(18 18)
\move(14 19.5) \rlvec(10 10)
\move(14 23.5) \rlvec(2 2)
\move(18 27.5) \rlvec(2 2)

\move(11 8) \textref h:R v:C \htext{\scriptsize{$0$}}
\move(15 8) \textref h:R v:C \htext{\scriptsize{$1$}}
\move(19 8) \textref h:R v:C \htext{\scriptsize{$2$}}
\move(23 8) \textref h:R v:C \htext{\scriptsize{$3$}}
\move(27 8) \textref h:R v:C \htext{\scriptsize{$4$}}
\move(31 8) \textref h:R v:C \htext{\scriptsize{$5$}}

\move(9 14) \textref h:R v:C \htext{\scriptsize{$1$}}
\move(9 18) \textref h:R v:C \htext{\scriptsize{$2$}}
\move(9 22) \textref h:R v:C \htext{\scriptsize{$3$}}
\move(9 26) \textref h:R v:C \htext{\scriptsize{$4$}}
\move(9 30) \textref h:R v:C \htext{\scriptsize{$5$}}

}
\caption{There are four tunnels in the \s path.}\label{tunnel}
\end{figure}
\end{ex}

The bijection $\phi$ can be easily described via the linear
representation of $\pi$. First in  $\gpi$,
add a mark right before each  singly covered minimal element
except $1$.
The bijection $\phi$
 transforms this marked linear representation of $\pi$ into a
lattice path by going through the vertices  from left to right, and replacing
each left end of an arc with an $N$-step, each right end of an arc with
an $E$-step, and each mark with a $D$-step.

\begin{ex}
Let
$\pi=\{\{1,6\},\{2,3,5\},\{3,4\}\}$. The marked linear representation is \\

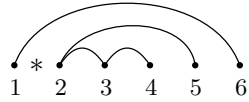
\begin{figure}[h,t]
\centertexdraw{
\drawdim mm

\linewd 0.2
\move(10 10) \clvec (14 21)(36 21)(40 10)
\move(16 10) \clvec (19 17)(33 17)(34 10)
\move(16 10) \clvec (18 13)(20 13)(22 10)
\move(22 10) \clvec (24 13)(26 13)(28 10)

\move(14 10) \textref h:R v:C \htext{$*$}

\move(10 10) \fcir f:0 r:0.4
\move(16 10) \fcir f:0 r:0.4
\move(22 10) \fcir f:0 r:0.4
\move(28 10) \fcir f:0 r:0.4
\move(34 10) \fcir f:0 r:0.4
\move(40 10) \fcir f:0 r:0.4

\move(11 7) \textref h:R v:C \htext{\footnotesize{$1$}}
\move(17 7) \textref h:R v:C \htext{\footnotesize{$2$}}
\move(23 7) \textref h:R v:C \htext{\footnotesize{$3$}}
\move(29 7) \textref h:R v:C \htext{\footnotesize{$4$}}
\move(35 7) \textref h:R v:C \htext{\footnotesize{$5$}}
\move(41 7) \textref h:R v:C \htext{\footnotesize{$6$}}
}
\end{figure}
The following steps yield the  corresponding \ls path.

\begin{figure}[h,t]
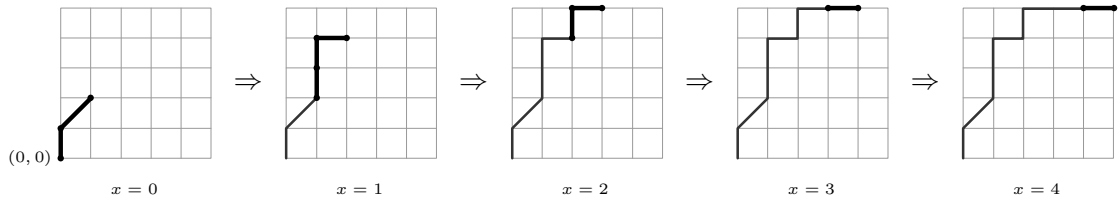

\centertexdraw{
\drawdim mm

\linewd 0.05 \setgray 0.6
\move(10 10) \rlvec(0 20)\rlvec(20 0)\rlvec(0 -20)\rlvec(-20 0)
\move(14 10) \rlvec(0 20)
\move(18 10) \rlvec(0 20)
\move(22 10) \rlvec(0 20)
\move(26 10) \rlvec(0 20)
\move(10 14) \rlvec(20 0)
\move(10 18) \rlvec(20 0)
\move(10 22) \rlvec(20 0)
\move(10 26) \rlvec(20 0)

\move(40 10) \rlvec(0 20)\rlvec(20 0)\rlvec(0 -20)\rlvec(-20 0)
\move(44 10) \rlvec(0 20)
\move(48 10) \rlvec(0 20)
\move(52 10) \rlvec(0 20)
\move(56 10) \rlvec(0 20)
\move(40 14) \rlvec(20 0)
\move(40 18) \rlvec(20 0)
\move(40 22) \rlvec(20 0)
\move(40 26) \rlvec(20 0)

\move(70 10) \rlvec(0 20)\rlvec(20 0)\rlvec(0 -20)\rlvec(-20 0)
\move(74 10) \rlvec(0 20)
\move(78 10) \rlvec(0 20)
\move(82 10) \rlvec(0 20)
\move(86 10) \rlvec(0 20)
\move(70 14) \rlvec(20 0)
\move(70 18) \rlvec(20 0)
\move(70 22) \rlvec(20 0)
\move(70 26) \rlvec(20 0)

\move(100 10) \rlvec(0 20)\rlvec(20 0)\rlvec(0 -20)\rlvec(-20 0)
\move(104 10) \rlvec(0 20)
\move(108 10) \rlvec(0 20)
\move(112 10) \rlvec(0 20)
\move(116 10) \rlvec(0 20)
\move(100 14) \rlvec(20 0)
\move(100 18) \rlvec(20 0)
\move(100 22) \rlvec(20 0)
\move(100 26) \rlvec(20 0)

\move(130 10) \rlvec(0 20)\rlvec(20 0)\rlvec(0 -20)\rlvec(-20 0)
\move(134 10) \rlvec(0 20)
\move(138 10) \rlvec(0 20)
\move(142 10) \rlvec(0 20)
\move(146 10) \rlvec(0 20)
\move(130 14) \rlvec(20 0)
\move(130 18) \rlvec(20 0)
\move(130 22) \rlvec(20 0)
\move(130 26) \rlvec(20 0)

\linewd 0.3 \setgray 0.2
\lpatt()
\move(40 10) \rlvec(0 4) \rlvec(4 4)

\move(70 10) \rlvec(0 4) \rlvec(4 4) \rlvec(0 8) \rlvec(4 0)

\move(100 10) \rlvec(0 4) \rlvec(4 4) \rlvec(0 8) \rlvec(4 0) \rlvec(0 4) \rlvec(4 0)

\move(130 10)\rlvec(0 4) \rlvec(4 4) \rlvec(0 8) \rlvec(4 0) \rlvec(0 4) \rlvec(8 0)

\linewd 0.6  \setgray 0
\move(10 10) \rlvec(0 4) \rlvec(4 4)

\move(44 18) \rlvec(0 8) \rlvec(4 0)

\move(82 30) \rlvec(-4 0) \rlvec(0 -4)

\move(112 30) \rlvec(4 0)

\move(146 30) \rlvec(4 0)

\move(10 10) \fcir f:0 r:0.4
\move(10 14) \fcir f:0 r:0.4
\move(14 18) \fcir f:0 r:0.4
\move(9 10) \textref h:R v:C \htext{\tiny{$(0,0)$}}

\move(44 18) \fcir f:0 r:0.4
\move(44 22) \fcir f:0 r:0.4
\move(44 26) \fcir f:0 r:0.4
\move(48 26) \fcir f:0 r:0.4

\move(78 26) \fcir f:0 r:0.4
\move(82 30) \fcir f:0 r:0.4
\move(78 30) \fcir f:0 r:0.4

\move(112 30) \fcir f:0 r:0.4
\move(116 30) \fcir f:0 r:0.4

\move(146 30) \fcir f:0 r:0.4
\move(150 30) \fcir f:0 r:0.4

\move(23 6) \textref h:R v:C \htext{\tiny{$x=0$}}
\move(53 6) \textref h:R v:C \htext{\tiny{$x=1$}}
\move(83 6) \textref h:R v:C \htext{\tiny{$x=2$}}
\move(113 6) \textref h:R v:C \htext{\tiny{$x=3$}}
\move(143 6) \textref h:R v:C \htext{\tiny{$x=4$}}

\move(37 20) \textref h:R v:C \htext{$\Rightarrow$}
\move(67 20) \textref h:R v:C \htext{$\Rightarrow$}
\move(97 20) \textref h:R v:C \htext{$\Rightarrow$}
\move(127 20) \textref h:R v:C \htext{$\Rightarrow$}

} \caption{The steps of $\phi$ that yield the
corresponding \s path.}
\end{figure}
\end{ex}

\begin{ex} The
elements of $NCL(3)$, their marked linear representations,  and the
corresponding \s paths.

\begin{figure}[ht]
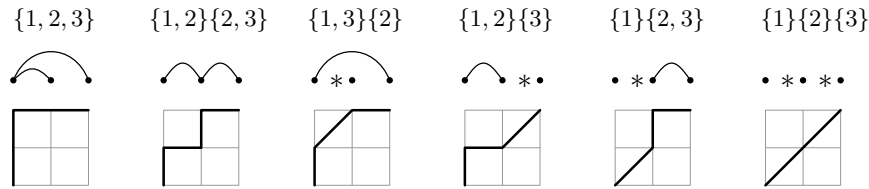

\centertexdraw{
\drawdim mm

\linewd 0.05 \setgray 0.6
\move(10 10) \rlvec(0 10)\rlvec(10 0)\rlvec(0 -10)\rlvec(-10 0)
\move(15 10) \rlvec(0 10)
\move(10 15) \rlvec(10 0)

\move(30 10) \rlvec(0 10)\rlvec(10 0)\rlvec(0 -10)\rlvec(-10 0)
\move(35 10) \rlvec(0 10)
\move(30 15) \rlvec(10 0)

\move(50 10) \rlvec(0 10)\rlvec(10 0)\rlvec(0 -10)\rlvec(-10 0)
\move(55 10) \rlvec(0 10)
\move(50 15) \rlvec(10 0)

\move(70 10) \rlvec(0 10)\rlvec(10 0)\rlvec(0 -10)\rlvec(-10 0)
\move(75 10) \rlvec(0 10)
\move(70 15) \rlvec(10 0)

\move(90 10) \rlvec(0 10)\rlvec(10 0)\rlvec(0 -10)\rlvec(-10 0)
\move(95 10) \rlvec(0 10)
\move(90 15) \rlvec(10 0)

\move(110 10) \rlvec(0 10)\rlvec(10 0)\rlvec(0 -10)\rlvec(-10 0)
\move(115 10) \rlvec(0 10)
\move(110 15) \rlvec(10 0)

\linewd 0.3 \setgray 0
\move (10 10) \rlvec(0 10) \rlvec(10 0)

\move (30 10) \rlvec(0 5) \rlvec(5 0) \rlvec(0 5) \rlvec(5 0)

\move (50 10) \rlvec(0 5) \rlvec(5 5) \rlvec(5 0)

\move (70 10) \rlvec(0 5) \rlvec(5 0) \rlvec(5 5)

\move (90 10) \rlvec(5 5) \rlvec(0 5) \rlvec(5 0)

\move (110 10) \rlvec(10 10)

\linewd 0.2 \setgray 0

\move(21 32) \textref h:R v:C \htext{\footnotesize{$\{1,2,3\}$}}
\move(44 32) \textref h:R v:C \htext{\footnotesize{$\{1,2\}\{2,3\}$}}
\move(62 32) \textref h:R v:C \htext{\footnotesize{$\{1,3\}\{2\}$}}
\move(82 32) \textref h:R v:C \htext{\footnotesize{$\{1,2\}\{3\}$}}
\move(102 32) \textref h:R v:C \htext{\footnotesize{$\{1\}\{2,3\}$}}
\move(124 32) \textref h:R v:C \htext{\footnotesize{$\{1\}\{2\}\{3\}$}}

\move(10 24) \clvec (12 29)(18 29)(20 24)
\move(10 24) \clvec (12 26)(13 26)(15 24)
\move(10 24) \fcir f:0 r:0.4
\move(15 24) \fcir f:0 r:0.4
\move(20 24) \fcir f:0 r:0.4

\move(30 24) \clvec (32 27)(33 27)(35 24)
\move(35 24) \clvec (37 27)(38 27)(40 24)
\move(30 24) \fcir f:0 r:0.4
\move(35 24) \fcir f:0 r:0.4
\move(40 24) \fcir f:0 r:0.4

\move(54 24) \textref h:R v:C \htext{$*$}
\move(50 24) \clvec (52 29)(58 29)(60 24)
\move(50 24) \fcir f:0 r:0.4
\move(55 24) \fcir f:0 r:0.4
\move(60 24) \fcir f:0 r:0.4

\move(79 24) \textref h:R v:C \htext{$*$}
\move(70 24) \clvec (72 27)(73 27)(75 24)
\move(70 24) \fcir f:0 r:0.4
\move(75 24) \fcir f:0 r:0.4
\move(80 24) \fcir f:0 r:0.4

\move(94 24) \textref h:R v:C \htext{$*$}
\move(95 24) \clvec (97 27)(98 27)(100 24)
\move(90 24) \fcir f:0 r:0.4
\move(95 24) \fcir f:0 r:0.4
\move(100 24) \fcir f:0 r:0.4

\move(114 24) \textref h:R v:C \htext{$*$}
\move(119 24) \textref h:R v:C \htext{$*$}
\move(110 24) \fcir f:0 r:0.4
\move(115 24) \fcir f:0 r:0.4
\move(120 24) \fcir f:0 r:0.4

}
\caption{The elements of $NCL(3)$ and their corresponding \s paths.}\label{f_NCL(3)}
\end{figure}
\end{ex}

A {\it peak} of the \s path is a pair of consecutive $NE$ steps and
a {\it valley} is a pair of consecutive $EN$ steps. The following
results are well-known for \s paths.

\begin{prop}\label{schroder}
\mbox{ } \par
\begin{enumerate}
\item Let $p(n,k)$ be the number of \s paths of length $n$ with $k$
   peaks. Then $ p(n,k)=C_{n-k}{2n-k\choose k}={n \choose k}{2n-k
   \choose n-1}/n$ where the $C_n$ is the $n$-th Catalan number
 $\frac{1}{n+1}{2n\choose n}$. It is also the number of \s paths
 with $k$ $D$-steps. Let $p(0,0)=1$. The generating function for $p(n,k)$ is
\begin{eqnarray} \label{peak}
\sum_{n, k \geq 0} p(n,k) x^nt^k = \frac{1-tx-\sqrt{(1-tx)^2-4x}}{2}.
\end{eqnarray}
\item Let $v(n,k)$ be the number of \s paths of length $n$  with k
valleys. It also counts the number of \s paths of length $n$
with $k$ $NN$-steps. Let $v(0,0)=1$.
The generating function $V(x, t)=\sum_{n,k \geq
  0} v(n,k) x^nt^k$ satisfies
$$
x(t+x-tx) V(x,t)^2-(1-2x+tx)V(x,t)+1=0.
$$
Explicitly,
\begin{eqnarray} \label{valley}
V(x,t)=\frac{-1+2x-tx+\sqrt{1-4x-2tx+t^2x^2}}{2(-tx-x^2+x^2t)}.
\end{eqnarray}
\item Let $d(n,k)$ be the number of \s paths of length $n$, containing
$k$ $D$'s not preceded by an $E$.
Let $d(0,0)=1$. The generating function
$D(x,t)=\sum_{n,k \geq 0} D(n,k)x^nt^k$ satisfies
$$
D(x,t)=1+txD(x,t)+x(1+x-tx)D(x,t)^2.
$$
Explicitly,
\begin{eqnarray} \label{singleton}
D(x,t)=\frac{1-tx-\sqrt{(1-tx)^2-4x(1+x-tx)}}{2x(1+x-tx)}.
\end{eqnarray}
\end{enumerate}
\end{prop}
These results can be found,
for example,  in OEIS \cite{Sloane}, Sequence A060693 for Statement 1,
A101282 for Statement 2, and A108916 for Statement 3.

The correspondence $\pi \rightarrow \phi(\pi)$ between noncrossing
linked partitions of $[n+1]$ and  \s paths of length $n$ allows us to
deduce a number of  properties for  noncrossing linked partitions.
It is easily seen that the properties  for an element of $[n+1]$ in a
noncrossing linked partition listed on the left  correspond to the given
steps of \s paths, listed on the right.

\begin{center}
\begin{tabular}{|c|c|} \hline
$NCL(n+1)$        & steps in \s paths \\ \hline
singly covered minimal element $i$,  $i \neq 1$ &   D \\ \hline
doubly covered element                 &   EN \\ \hline
 singleton block $\{i\}$, $i \neq 1$          &  D not followed by an N \\ \hline
$i \in B$ where $\min(B)=i-1$          &   NE \\  \hline
\end{tabular}
\end{center}

From  Prop. \ref{schroder} and the obvious symmetry between steps $ED$
and $DN$ there follows:
\begin{prop}
\begin{enumerate}
\item The number of noncrossing linked partitions of $[n+1]$ with $k$
   singly covered minimal elements $i$ where $i \neq 1$ is equal to
   $p(n,k)=C_{n-k}{2n-k\choose k}={n \choose k}{2n-k  \choose n-1}/n$. It is
  also counts the number of noncrossing linked partitions on $[n+1]$
 with $k$ elements $x$ such that $x, x-1$ lie in a block $B$ with
 $x-1=\min(B)$. The generating function of $p(n,k)$ is given by
   Eqn. \eqref{peak}.
\item The number of noncrossing linked partitions of $[n+1]$ with $k$
   doubly covered elements is $v(n,k)$, whose generating function is
   given by Eqn. \eqref{valley}.
\item The number of noncrossing linked partitions on $[n+1]$ with $k$
 singleton blocks $\{i\}$ where $i \neq 1$ is $d(n,k)$,
whose generating function is
given by Eqn. \eqref{singleton}.
\end{enumerate}
\end{prop}

At the end of this section,
we count the noncrossing linked partitions by the number of
blocks, using the recurrence \eqref{rec-for}.

\begin{prop}
Let $b(n,k)$ be the number of noncrossing linked partition of $[n]$
with $k$ blocks. Let $B(x,t)=1+\sum_{n, k \geq 1} b(n,k)x^n t^k$.
Then $B(x,t)$ satisfies the equation
\begin{eqnarray} \label{eqn_bnk}
(1+x-tx)B(x,t)^2+(2tx-x-3)B(x,t)+2=0.
\end{eqnarray}
Explicitly,
\begin{eqnarray} \label{exact-for}
B(x,t)=\frac{3+x-tx-\sqrt{(2tx-x-3)^2-8(1+x-tx)}}{2(1+x-tx)}.
\end{eqnarray}
\end{prop}
\begin{proof}
We derive a recurrence for $b(n,k)$. Given $\pi \in NCL(n)$ with $k$
blocks,
again let $i=\min(B[n])$.
First, there are $b(n-1, k-1)$ many noncrossing linked
partitions in $NCL(n)$ such that $i=n$.
Otherwise, assume $1 \leq i \leq n-1$. As in the proof of
Prop. \ref{recurrence}, $\pi$ is a union of two noncrossing linked
partitions, $\pi_1$ of $[i]$, and $\pi_2$ of $\{i, \dots, n\}$ with
$i$ and $n$ lying in the same block. Assume $\pi_1$ has $t$ blocks,
and $\pi_2$ has $r$ blocks. If $i$ is a singleton of $\pi_1$, then
$\pi$ has $t-1+r$ many blocks; if $i$ is not a singleton of $\pi_1$,
then $\pi$ has $t+r$ many blocks. Finally, note that $\pi_2$ can be
obtained by taking any noncrossing linked partition on $\{i, i+1,
\dots, n-1\}$, and then adding $n$ to the block containing
$i$. Combining the above, we get the recurrence
\begin{eqnarray} \label{rec-b_nk}
b(n,k)&=& b(n-1, k-1) \nonumber \\
       &  & +\sum_{i=1}^{n-1} \sum_{r+t=k} b(n-i, r)\Big[b(i-1,
t)+b(i,t)-b(i-1, t-1)\Big].
\end{eqnarray}
Note $b(1,1)=b(2,1)=b(2,2)=1$, and $b(n,0)=0$ for all $n\geq 1$.
If we set $b(0,0)=1$, then the recurrence \eqref{rec-b_nk} holds for
all $n \geq 1$ and $k \geq 1$.
Now multiply both sides of \eqref{rec-b_nk} by $x^n t^k$, and sum
over all $n, k$.  Noticing that for any sequence $g_i$,
$$
 \sum_{n \geq 1} \sum_{i=1}^{n-1} g_ig_{n-i} x^n = \left(\sum_{n \geq 0}
 g_nx^n \right)^2-2g_0 \sum_{n \geq 0} g_nx^n,
$$
we get the equation \eqref{eqn_bnk} and hence the formula
\eqref{exact-for}.
\end{proof}


\section{Linked partitions}
In this section we study linked partitions. Recall that
a  linked partition of $[n]$ is a collection of pairwise nearly
disjoint subsets whose union is $[n]$.
The set of all the linked partitions on $[n]$ is denoted by $LP(n)$,
whose cardinality is  $lp_n$.

It is not hard to see that $lp_n=n!$. Instead of merely giving a
counting argument,  we present a one-to-one correspondence between
the set of linked partitions of $[n]$ and the set of increasing trees
on $n+1$ labeled vertices. The latter is a geometric representation
for permutations, originally developed by the French, and
outlined in the famous textbook \cite[Chap.1.3]{Stan}.
Many properties of linked partitions can be trivially deduced from
this correspondence. As a sample, we list the results involving
the signless Stirling numbers  and the Eulerian
numbers. At the end of this section we  give the joint distribution
for two statistics,
2-crossings and 2-nestings, over linked partitions with given
sets of lefthand and righthand endpoints.

\begin{defi}
An {\it increasing tree} on $n+1$ labeled vertices is a rooted tree on
vertices $0,1,\ldots,n$ such that for any vertex $i$, $i < j$
if $i$ is a successor of $j$.
\end{defi}

\begin{thm}
There is a one-to-one correspondence between the set of linked
partitions of $[n]$ and the set of increasing trees on $n+1$ labeled vertices.
\end{thm}
\begin{proof}
We use the  linear representation $\gpi$ for linked partitions.
Recall that  for $\pi \in LP(n)$,
$\gpi$ is the graph with $n$ dots listed in a horizontal line with
labels $1, 2, \dots, n$, where  $i$ and $j$ are connected by
 an arc if and only of
$j$ lies in a block $B$ with $i=\min(B)$.
To get an increasing tree, one simply adds a root 0
to $\gpi$
which connects to all the singly covered minimal elements of $\pi$.
This defines
an increasing tree on $[n] \cup \{0\}$, where the children of root $0$
are
those singly covered
minimal elements, and $j$ is a child of $i$ if and only if
$(i,j)$ is an arc of $\gpi$ and $j>i$.
\end{proof}

\begin{ex}\label{Ex_lP}
Let $\pi=\{\{126\}, \{248\}, \{3\}, \{57\}\}$.  The
singly covered minimal elements are $\{1, 3, 5\}$. The
corresponding increasing tree is given in Figure \ref{LP-Itree}.

\begin{figure}[h,t]
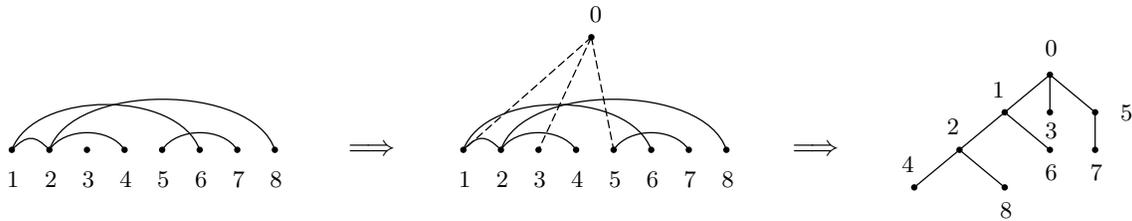

\centertexdraw{ \drawdim mm

\linewd 0.2
\move(10 15) \clvec (12 17)(13 17)(15 15) 
\move(10 15) \clvec (13 23)(32 23)(35 15) 
\move(15 15) \clvec (17 18)(23 18)(25 15) 
\move(15 15) \clvec (18 24)(42 24)(45 15) 
\move(30 15) \clvec (32 18)(38 18)(40 15) 

\move(10 15) \fcir f:0 r:0.4
\move(15 15) \fcir f:0 r:0.4
\move(20 15) \fcir f:0 r:0.4
\move(25 15) \fcir f:0 r:0.4
\move(30 15) \fcir f:0 r:0.4
\move(35 15) \fcir f:0 r:0.4
\move(40 15) \fcir f:0 r:0.4
\move(45 15) \fcir f:0 r:0.4

\move(11 11) \textref h:R v:C \htext{\footnotesize{1}}
\move(16 11) \textref h:R v:C \htext{\footnotesize{2}}
\move(21 11) \textref h:R v:C \htext{\footnotesize{3}}
\move(26 11) \textref h:R v:C \htext{\footnotesize{4}}
\move(31 11) \textref h:R v:C \htext{\footnotesize{5}}
\move(36 11) \textref h:R v:C \htext{\footnotesize{6}}
\move(41 11) \textref h:R v:C \htext{\footnotesize{7}}
\move(46 11) \textref h:R v:C \htext{\footnotesize{8}}

\move(61 15) \textref h:R v:C \htext{$\Longrightarrow$}

\move(70 15) \clvec (72 17)(73 17)(75 15) 
\move(70 15) \clvec (73 23)(92 23)(95 15) 
\move(75 15) \clvec (77 18)(83 18)(85 15) 
\move(75 15) \clvec (78 24)(102 24)(105 15) 
\move(90 15) \clvec (92 18)(98 18)(100 15) 

\move(70 15) \fcir f:0 r:0.4
\move(75 15) \fcir f:0 r:0.4
\move(80 15) \fcir f:0 r:0.4
\move(85 15) \fcir f:0 r:0.4
\move(90 15) \fcir f:0 r:0.4
\move(95 15) \fcir f:0 r:0.4
\move(100 15) \fcir f:0 r:0.4
\move(105 15) \fcir f:0 r:0.4

\move(71 11) \textref h:R v:C \htext{\footnotesize{1}}
\move(76 11) \textref h:R v:C \htext{\footnotesize{2}}
\move(81 11) \textref h:R v:C \htext{\footnotesize{3}}
\move(86 11) \textref h:R v:C \htext{\footnotesize{4}}
\move(91 11) \textref h:R v:C \htext{\footnotesize{5}}
\move(96 11) \textref h:R v:C \htext{\footnotesize{6}}
\move(101 11) \textref h:R v:C \htext{\footnotesize{7}}
\move(106 11) \textref h:R v:C \htext{\footnotesize{8}}

\move(87 30) \fcir f:0 r:0.4
\lpatt(1 0.6)
\move(87 30) \rlvec(-17 -15)
\move(87 30) \rlvec(-7 -15)
\move(87 30) \rlvec(3 -15)

\move(88.5 33) \textref h:R v:C \htext{\footnotesize{0}}


\move(120 15) \textref h:R v:C \htext{$\Longrightarrow$}

\lpatt()
\move(130 10) \rlvec(18 15) \rlvec(6 -5)
\move(136 15) \rlvec(6 -5)
\move(142 20) \rlvec(6 -5)
\move(148 25) \rlvec(0 -5)
\move(154 20) \rlvec(0 -5)

\move(130 10) \fcir f:0 r:0.4
\move(136 15) \fcir f:0 r:0.4
\move(142 10) \fcir f:0 r:0.4
\move(142 20) \fcir f:0 r:0.4
\move(148 15) \fcir f:0 r:0.4
\move(148 20) \fcir f:0 r:0.4
\move(148 25) \fcir f:0 r:0.4
\move(154 20) \fcir f:0 r:0.4
\move(154 15) \fcir f:0 r:0.4

\move(130 13) \textref h:R v:C \htext{\footnotesize{$4$}}
\move(136 18) \textref h:R v:C \htext{\footnotesize{$2$}}
\move(142 23) \textref h:R v:C \htext{\footnotesize{$1$}}
\move(149 28.5) \textref h:R v:C \htext{\footnotesize{$0$}}
\move(143 7) \textref h:R v:C \htext{\footnotesize{$8$}}
\move(149 12) \textref h:R v:C \htext{\footnotesize{$6$}}
\move(149 17.5) \textref h:R v:C \htext{\footnotesize{$3$}}
\move(159 20) \textref h:R v:C \htext{\footnotesize{$5$}}
\move(155 12) \textref h:R v:C \htext{\footnotesize{$7$}}

}
\caption{A linked partition and the corresponding increasing tree}
\label{LP-Itree}
\end{figure}
\end{ex}

The following properties of  increasing trees are listed in
Proposition 1.3.16 of
\cite{Stan}.

\begin{prop}
\begin{enumerate}
\item The number of increasing trees on $n+1$ labeled vertices is $n!$.
\item The number of such trees for which the root has $k$ successors
is the signless Stirling number $c(n,k)$ (of the first kind).
\item The number of such trees with $k$ endpoints is the Eulerian
number $A(n,k)$.
\end{enumerate}
\end{prop}

Let $\beta(\pi)=|\{i \text{ is singly covered and $i \neq
  \min(B[i])$}\}| +  |\{i:  \text{  singleton block of $\pi$ } \}|$.
\begin{cor}
\begin{enumerate}
\item The number of linked partitions of $[n]$ is $n!$.

\item The number of  linked partitions of $[n]$ with $k$ singly covered
minimal elements is the signless Stirling number $c(n,k)$ (of the
first kind).

\item The number of linked partitions of $[n]$ with $\beta(\pi)=k$ is the
Eulerian number $A(n,k)$.
\end{enumerate}
\end{cor}

For a linked partition $\pi$ with linear representation $\gpi$,
we say that
two arcs $(i_1,j_1)$ and $(i_2,j_2)$  form a
{\it 2-crossing} if $i_1<i_2<j_1<j_2$; they form a {\it 2-nesting} if
$i_1<i_2<j_2<j_1$. Denoted   by  $cr_2(\pi)$ and $ne_2(\pi)$ the
number of 2-crossings and 2-nestings of $\pi$,  respectively. For example,
the linked partition $\pi=\{\{126\}, \{248\}, \{3\}, \{57\}\}$ in
Example \ref{Ex_lP} has $cr_2(\pi)=2$ and $ne_2(\pi)=2$.

Given $\pi\in LP(n)$, define two multiple sets
\begin{eqnarray*}
\text{left}(\pi)~=~\{\text{lefthand endpoints of arcs of }\pi\},\\
\text{right}(\pi)~=~\{\text{righthand endpoints of arcs of }\pi\}.
\end{eqnarray*}
For example, for $\pi=\{\{126\}, \{248\}, \{3\}, \{57\}\}$,
left$(\pi)=\{1,1,2,2,5\}$, and right$(\pi)=\{2,6,4,7,8\}$.
Clearly, each element of $\text{right}(\pi)$ has multiplicity $1$.

Fix $S$ and $T$ where $S$ is a multi-subset of $[n]$,
$T$ is a subset of $[n]$, and  $|S|=|T|$. Let $LP_n(S,T)$ be the set
\{$\pi\in LP(n):$ left$(\pi)=S$, right$(\pi)=T$\}.
We prove that over each set $LP_n(S, T)$, the statistics
 $cr_2(\pi)$ and $ne_2(\pi)$ have a symmetric
joint distribution. Explicitly, let
 $S=\{a_1^{r_1}, a_2^{r_2}, \ldots, a_m^{r_m}\}$
with $a_1 < a_2 < \cdots < a_m$.
For each $1 \leq i \leq m$, let
$h(i)=|\{j \in T: j >  a_i\}|-|\{j \in S: j >a_i\}|$.
We have

\begin{thm}\label{thm_cr_ne}
\begin{equation}\label{LP-cr-ne}
\sum_{\pi\in LP_n(S,T)}x^{cr_2(\pi)}y^{ne_2(\pi)}
= \sum_{\pi\in LP_n(S,T)}x^{ne_2(\pi)}y^{cr_2(\pi)}
= \prod_{i=1}^m y^{r_i h(i)-r_i^2} \left.\mathbf{ {h(i) \choose r_i }}\right|_{q=x/y},
\end{equation}
where $\mathbf{{n \choose m}}$ is the q-binomial coefficient
$$
\mathbf{ {n \choose m}}=\frac{\mathbf{(n)!}}{\mathbf{ (m)!(n-m)!}}
=\frac{ (q^n-1)(q^n-q)\cdots (q^n-q^{m-1})}{(q^m-1)(q^m-q)\cdots
  (q^m-q^{m-1})}.$$
\end{thm}
\begin{proof}
For $1 \leq m \leq n$, denote by ${ [n] \choose m}$ the set of integer
sequences
$(x_1, x_2, \dots, x_m)$ such that $1 \leq x_1 < x_2 < \dots < x_m \leq n$.
We give a bijection from the set of linked partitions in $LP_n(S, T)$
to the set $\prod_{i=1}^m {[h(i)] \choose r_i}$.

Given an element $\mathbf{s=(s_1, s_2, \dots, s_m)}$
 in $\prod_{i=1}^m {[h(i)] \choose r_i}$ where $\mathbf{s_i} \in
 {[h(i)] \choose r_i}$,  we construct a
linked partition $\pi$ by matching each lefthand endpoint in $S$ to a
righthand endpoint in $T$.
First, there are $r_m$ lefthand endpoints
at node $a_m$, and on its right there are $h(m)$ many righthand endpoints.
Assume $\mathbf{s_m}=(x_1, x_2, \dots, x_{r_m})$. We connect
the $r_m$ lefthand endpoints at node $a_m$ to the $x_1$-th, $x_2$-th, ...,
$x_{r_m}$-th righthand endpoints after node $a_m$.

In general, after matching lefthand endpoints at
nodes $a_{i+1}, \dots,
a_m$ to some righthand end-points, we process the $r_i$ lefthand endpoints at
node $a_i$.
At this stage there are exactly $h(i)$ many righthand endpoints available
after the node $a_i$.  List them by $1, 2, \dots, h(i)$, and match
the lefthand endpoints at $a_i$ to   the $y_1$-th, $y_2$-th, ...,
$y_{r_i}$-th  of them, if
$\mathbf{s_i}=(y_1, y_2, \dots, y_{r_i})$.

Continue the above procedure until each  lefthand endpoint is connected to
some righthand endpoint on its  right. This gives the desired bijection between
$LP_n(S, T)$ and $\prod_{i=1}^m {[h(i)] \choose r_i}$. In particular,
$LP_n(S, T)$ is nonempty if and only if for every
$a_i \in S$, $h(i) \geq r_i$.

\begin{ex}
Let $S=\{1,1,2,2,2, 3, 3\}$, $T=\{4,5,6,7,8,9,10\}$. Then
$h(1)=2$, $h(2)=5$, $h(3)=7$.  Figure \ref{LP(S,T)} illustrates how to
construct the linked partition for $\mathbf{s=(s_1,s_2,s_3)}$ where
$\mathbf{s_1}=(1,2),
\mathbf{s_2}=(2,4,5)$, and $\mathbf{s_3}=(3,6)$.

\begin{figure}[h,t]
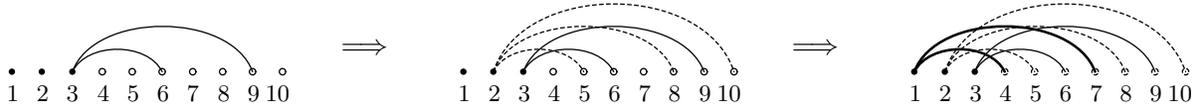

\centertexdraw{ \drawdim mm

\move(10 23) \textref h:R v:C \htext{}

\linewd 0.2
\move(10 10) \fcir f:0 r:0.4
\move(14 10) \fcir f:0 r:0.4
\move(18 10) \fcir f:0 r:0.4
\move(22 10) \lcir r:0.4
\move(26 10) \lcir r:0.4
\move(30 10) \lcir r:0.4
\move(34 10) \lcir r:0.4
\move(38 10) \lcir r:0.4
\move(42 10) \lcir r:0.4
\move(46 10) \lcir r:0.4

\move(11 7) \textref h:R v:C \htext{\footnotesize{1}}
\move(15 7) \textref h:R v:C \htext{\footnotesize{2}}
\move(19 7) \textref h:R v:C \htext{\footnotesize{3}}
\move(23 7) \textref h:R v:C \htext{\footnotesize{4}}
\move(27 7) \textref h:R v:C \htext{\footnotesize{5}}
\move(31 7) \textref h:R v:C \htext{\footnotesize{6}}
\move(35 7) \textref h:R v:C \htext{\footnotesize{7}}
\move(39 7) \textref h:R v:C \htext{\footnotesize{8}}
\move(43 7) \textref h:R v:C \htext{\footnotesize{9}}
\move(47 7) \textref h:R v:C \htext{\footnotesize{10}}

\move(18 10) \clvec (20 14)(28 14)(30 10) 
\move(18 10) \clvec (21 18)(39 18)(42 10) 


\move(60 13) \textref h:R v:C \htext{$\Longrightarrow$}

\move(70 10) \fcir f:0 r:0.4
\move(74 10) \fcir f:0 r:0.4
\move(78 10) \fcir f:0 r:0.4
\move(82 10) \lcir r:0.4
\move(86 10) \lcir r:0.4
\move(90 10) \lcir r:0.4
\move(94 10) \lcir r:0.4
\move(98 10) \lcir r:0.4
\move(102 10) \lcir r:0.4
\move(106 10) \lcir r:0.4

\move(71 7) \textref h:R v:C \htext{\footnotesize{1}}
\move(75 7) \textref h:R v:C \htext{\footnotesize{2}}
\move(79 7) \textref h:R v:C \htext{\footnotesize{3}}
\move(83 7) \textref h:R v:C \htext{\footnotesize{4}}
\move(87 7) \textref h:R v:C \htext{\footnotesize{5}}
\move(91 7) \textref h:R v:C \htext{\footnotesize{6}}
\move(95 7) \textref h:R v:C \htext{\footnotesize{7}}
\move(99 7) \textref h:R v:C \htext{\footnotesize{8}}
\move(103 7) \textref h:R v:C \htext{\footnotesize{9}}
\move(107 7) \textref h:R v:C \htext{\footnotesize{10}}

\move(78 10) \clvec (80 14)(88 14)(90 10) 
\move(78 10) \clvec (81 18)(99 18)(102 10) 

\lpatt(0.5 0.4)
\move(74 10) \clvec (76 14)(84 14)(86 10) 
\move(74 10) \clvec (77 18)(95 18)(98 10) 
\move(74 10) \clvec (78 22)(102 22)(106 10) 


\move(120 13) \textref h:R v:C \htext{$\Longrightarrow$}

\move(130 10) \fcir f:0 r:0.4
\move(134 10) \fcir f:0 r:0.4
\move(138 10) \fcir f:0 r:0.4
\move(142 10) \lcir r:0.4
\move(146 10) \lcir r:0.4
\move(150 10) \lcir r:0.4
\move(154 10) \lcir r:0.4
\move(158 10) \lcir r:0.4
\move(162 10) \lcir r:0.4
\move(166 10) \lcir r:0.4

\move(131 7) \textref h:R v:C \htext{\footnotesize{1}}
\move(135 7) \textref h:R v:C \htext{\footnotesize{2}}
\move(139 7) \textref h:R v:C \htext{\footnotesize{3}}
\move(143 7) \textref h:R v:C \htext{\footnotesize{4}}
\move(147 7) \textref h:R v:C \htext{\footnotesize{5}}
\move(151 7) \textref h:R v:C \htext{\footnotesize{6}}
\move(155 7) \textref h:R v:C \htext{\footnotesize{7}}
\move(159 7) \textref h:R v:C \htext{\footnotesize{8}}
\move(163 7) \textref h:R v:C \htext{\footnotesize{9}}
\move(167 7) \textref h:R v:C \htext{\footnotesize{10}}

\lpatt()
\move(138 10) \clvec (140 14)(148 14)(150 10) 
\move(138 10) \clvec (141 18)(159 18)(162 10) 

\lpatt(0.5 0.4)
\move(134 10) \clvec (136 14)(144 14)(146 10) 
\move(134 10) \clvec (137 18)(155 18)(158 10) 
\move(134 10) \clvec (138 22)(162 22)(166 10) 

\linewd 0.3
\lpatt()
\move(130 10) \clvec (132 14)(140 14)(142 10) 
\move(130 10) \clvec (133 18)(151 18)(154 10) 

} \caption{An illustration of the bijection from $LP_n(S, T)$ to
$\prod_{i=1}^m {[h(i)] \choose r_i}$.}
\label{LP(S,T)}
\end{figure}
\end{ex}

The numbers of 2-crossings and 2-nestings are easily expressed in
terms of  the
sequence $\mathbf{s}=\mathbf{(s_1, s_2, \dots, s_m)}$.
Assume $\mathbf{s_i}=(x_1, \dots, x_{r_i})$.  By the above construction,
the number of 2-crossings formed by arcs $jk$ and $a_i t$ with
$j < a_i < k < t$ is $\sum_{t=1}^{r_i} (x_t-t)$, and the number of
2-crossings formed by arcs $jk$ and $a_i t$ with $j < a_i < t < k$
is $\sum_{t=1}^{r_i} (h(i)-x_t-(r_i-t))$.
Hence $\mathbf{s_i}$ contributes a factor
$$
x^{\left(\sum_{t=1}^{r_i} x_t\right) -{r_i+1 \choose 2}}
y^{r_i h(i)-{r_i \choose 2}-\sum_{t=1}^{r_i} x_t}
$$
to the generating function
$\sum_{\pi\in LP_n(S,T)}x^{cr_2(\pi)}y^{ne_2(\pi)}$.
Since
$$
\sum_{(x_1, \dots, x_{r_i}) \in {[h(i)]\choose r_i}} q^{\left(\sum_t x_t\right) -{r_i +1 \choose 2}} =  \mathbf{ {h(i) \choose r_i }},
$$
and $\mathbf{s_i}$ are mutually independent, we have
$$
\sum_{\pi\in LP_n(S,T)}x^{cr_2(\pi)}y^{ne_2(\pi)}
= \left.\prod_{i=1}^m y^{r_i h(i)-r_i^2} \mathbf{ {h(i) \choose r_i }}\right|_{q=x/y}.
$$
The symmetry between $cr_2(\pi)$ and $ne_2(\pi)$ is obtained by
the involution $\tau: \mathbf{(s_1, s_2, \dots, s_m)}
 \rightarrow \mathbf{(t_1, t_2, \dots, t_m)}$ on
 $\prod_{i=1}^m {[h(i)] \choose r_i}$, where
 $\mathbf{t_i}=(h(i)+1-x_{r_i}, \dots, h(i)+1-x_1)$
if $\mathbf{s_i}=(x_1, \dots, x_{r_i})$.
\end{proof}

%
%

\section{Linked cycles}
As with matchings and
 set partitions,
one can define the \emph{crossing number} and the \emph{nesting
 number}  for a
given linked
partition. Unfortunately, these two statistics do not have the same
distribution
over all linked partitions of $[n]$.
Motivated by the work in \cite{CDDSY},
we want to find a suitable structure over which
the crossing number and nesting number have a symmetric joint
 distribution.
For this purpose
we introduce the notion of \emph{linked cycles}, which are   linked
 partitions
equipped with
a cycle structure on each of its block.
It turns out that the set of
 linked
cycles possesses
many interesting combinatorial properties.

\subsection{Two representations for linked cycles}

\begin{defi}
A {\it linked cycle} $\hat \pi$ on $[n]$ is a linked partition $\pi=\{B_1, \dots, B_k\}$
of $[n]$ where for each block $B_i$ the elements are arranged in a
cycle.
\end{defi}

We call each such block $B_i$ with the cyclic arrangement a {\it cycle} of $\hat \pi$.
The set of all linked cycles on $[n]$ is denoted by $LC(n)$.

For a set $B=\{b_1,b_2,\ldots,b_k\}$, we represent by
$(b_1,b_2,\ldots,b_k)$ the cycle $b_1 \goto b_2\goto \cdots\goto b_k
\goto b_1$. In writing a linked cycle $\hat \pi$, we use the
convention that: (a) each cycle of $\hat \pi$  is written with its
minimal element first, (b) the cycles are listed in increasing order
of their minimal elements. For example, for the linked partition
$\pi=\{ \{126\},\{248\}, \{3\}, \{57\}\}$ with cyclic orders $(1
\goto 2 \goto 6 \goto 1), (2 \goto 8 \goto 4 \goto 2), (3 \goto 3)$,
and $(5 \goto 7 \goto 5)$, the linked cycle $\hat \pi$ is written as
$\hat \pi=(126)(284)(3)(57)$.

Again the linked cycles may be represented by certain graphs.  Here we introduce
 two such graphical  representations.

\noindent {\bf 1. Cycle representation $\ghpi^c$}. \
Let $\hat \pi \in LC(n)$, the cycle representation
$\ghpi^c$ of $\hat \pi$ is  a directed graph  on
$[n]$ with arcs $(i,j)$ whenever $i$ and $j$ are  consecutive
elements in a cycle of $\hat \pi$.

In drawing the figures, we put elements of a cycle $C_i$ in a circle in
clockwise order.
If a cycle $C_i$ contains the minimal
element of a connected component of $\ghpi^c$, then we say that $C_i$ is the {\it root cycle} of
that component.

\begin{ex}
The cycle representation for the linked cycle $\pi=(126)(284)(3)(57)$
is given in Figure \ref{LC}. There are three connected
components and $(126)$ is the root cycle of the connected component
$(126)(284)$.

\begin{figure}[h,t]
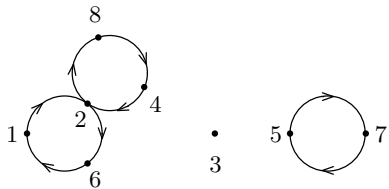

\centertexdraw{ \drawdim mm

\linewd 0.2

\move(10 10) \lcir r:5

\move(5 10) \fcir f:0 r:0.4
\move(13 14) \fcir f:0 r:0.4
\move(13 6) \fcir f:0 r:0.4

\move(4 10) \textref h:R v:C \htext{\footnotesize{$1$}}
\move(13 12) \textref h:R v:C \htext{\footnotesize{$2$}}
\move(15 4) \textref h:R v:C \htext{\footnotesize{$6$}}

\arrowheadsize l:1.4 w:1
\move(6 13) \arrowheadtype t:v \avec(7 14)
\move(15 10.5) \arrowheadtype t:v \avec(15 9.3)
\move(8 5.5) \arrowheadtype t:v \avec(7 6)

\move(16 18) \lcir r:5

\move(14.5 22.8) \fcir f:0 r:0.4
\move(20.6 16.2) \fcir f:0 r:0.4

\move(15 26) \textref h:R v:C \htext{\footnotesize{$8$}}
\move(23 14) \textref h:R v:C \htext{\footnotesize{$4$}}

\arrowheadsize l:1.4 w:1
\move(11 18) \arrowheadtype t:v \avec(11.2 19.2)
\move(20.4 20.5) \arrowheadtype t:v \avec(20.8 19.2)
\move(18.5 13.6) \arrowheadtype t:v \avec(17 13)

\move(30 10) \fcir f:0 r:0.4

\move(31 6) \textref h:R v:C \htext{\footnotesize{$3$}}

\move(45 10) \lcir r:5

\move(40 10) \fcir f:0 r:0.4
\move(50 10) \fcir f:0 r:0.4

\move(39 10) \textref h:R v:C \htext{\footnotesize{$5$}}
\move(53 10) \textref h:R v:C \htext{\footnotesize{$7$}}

\arrowheadsize l:1.4 w:1
\move(44.3 14.9) \arrowheadtype t:v \avec(45.7 14.9)
\move(45.7 5.1) \arrowheadtype t:v \avec(44.3 5.1)

} \caption{The cycle representation of the linked cycle
$\pi=(126)(248)(3)(57)$.}\label{LC}
\end{figure}
\end{ex}

\noindent {\bf 2. Linear representation $\ghpi^l$}. Let $\hat \pi\in
LC(n)$, the linear representation $\ghpi^l$ of $\hat \pi$ is  a
graph whose vertices lie on a horizontal line, and each vertex is of
one of the following kind:  (i) a lefthand endpoint, (ii) a
righthand endpoint, or (iii) an isolated point.

Start with $n$ vertices in a horizontal line labeled $1, 2, \dots, n$.
We define $\ghpi^l$ by first splitting vertices  as follows:
\begin{enumerate}
\item
 If $i$ is a singly covered minimal element of a cycle with $k+1$
 elements, and $k \geq 0$, split the  vertex $i$ into $k$ vertices labeled
$i^{(1)}, i^{(2)}, \ldots, i^{(k)}$;
\item If $i$ is  singly covered, but not a minimal element,
 replace the vertex $i$  by a vertex labeled $i^{(0)}$;
\item If $i$ is  doubly covered, and is the minimal element of a cycle
  of size $k+1$,  split the vertex $i$  into
$k+1$  vertices labeled $i^{(0)}, i^{(1)}, i^{(2)}, \ldots, i^{(k)}$;
\end{enumerate}
For a cycle $C_i=(i_1 i_2\ldots i_{t_i})$,
if $t_i \geq 2$, then we have created vertices with labels
$i_1^{(1)}, i_1^{(2)}, \ldots, i_1^{(t_i-1)}$.  Add arcs
$(i_1^{(1)}, i_{t_i}^{(0)})$,
$(i_1^{(2)},i_{t_i-1}^{(0)}),~\ldots, (i_1^{(t_i-1)},i_2^{(0)})$.
Do this for each cycle of $\hp$, and the resulting graph is $\gpi^l$.

\begin{ex}
The linear representation for the linked cycle $\hat \pi=(126)(284)(3)(57)$.
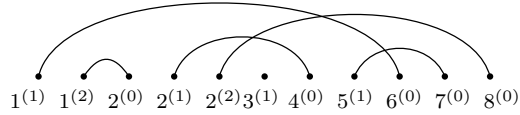
\vspace{.8cm}

\begin{figure}[h,t]
\centertexdraw{ \drawdim mm

\linewd 0.2

\move(104 10) \clvec (107 23)(149 23)(152 10)
\move(110 10) \clvec (112 13)(114 13)(116 10)
\move(122 10) \clvec (124 17)(138 17)(140 10)
\move(128 10) \clvec (131 21)(161 21)(164 10)
\move(146 10) \clvec (148 15)(156 15)(158 10)

\move(104 10) \fcir f:0 r:0.4
\move(110 10) \fcir f:0 r:0.4
\move(116 10) \fcir f:0 r:0.4
\move(122 10) \fcir f:0 r:0.4
\move(128 10) \fcir f:0 r:0.4
\move(134 10) \fcir f:0 r:0.4
\move(140 10) \fcir f:0 r:0.4
\move(146 10) \fcir f:0 r:0.4
\move(152 10) \fcir f:0 r:0.4
\move(158 10) \fcir f:0 r:0.4
\move(164 10) \fcir f:0 r:0.4

\move(105 7) \textref h:R v:C \htext{\scriptsize{$1^{(1)}$}}
\move(111.5 7) \textref h:R v:C \htext{\scriptsize{$1^{(2)}$}}
\move(118 7) \textref h:R v:C \htext{\scriptsize{$2^{(0)}$}}
\move(124.5 7) \textref h:R v:C \htext{\scriptsize{$2^{(1)}$}}
\move(131 7) \textref h:R v:C \htext{\scriptsize{$2^{(2)}$}}
\move(136 7) \textref h:R v:C \htext{\scriptsize{$3^{(1)}$}}
\move(142 7) \textref h:R v:C \htext{\scriptsize{$4^{(0)}$}}
\move(148.5 7) \textref h:R v:C \htext{\scriptsize{$5^{(1)}$}}
\move(155 7) \textref h:R v:C \htext{\scriptsize{$6^{(0)}$}}
\move(161.5 7) \textref h:R v:C \htext{\scriptsize{$7^{(0)}$}}
\move(168 7) \textref h:R v:C \htext{\scriptsize{$8^{(0)}$}}
}
\label{linear rep.}
\end{figure}
\end{ex}

As for the linked partitions,  sometimes it is useful to distinguish the
set of singly covered minimal elements of  $\hat \pi$. A vertex $i$ is
a singly covered minimal element of $\hat \pi$ if and only if
in the linear representation  $\ghpi^l$, there is no vertex $i^{(0)}$.
The {\it marked linear representation} of $\hat \pi$ is obtained
from $\ghpi^l$ by adding a mark before $i^{(1)}$  for each
singly covered minimal $i$, except for $i=1$. For example, for
$\hp$ in Example \ref{linear rep.} we should add marks before vertices
$3^{(1)}$ and $5^{(1)}$.

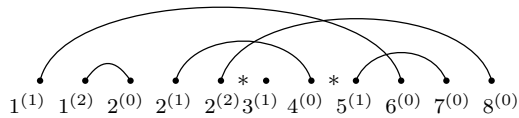
\begin{figure}[h,t]
\centertexdraw{ \drawdim mm

\linewd 0.2
\move(100 25) \textref h:R v:C \htext{}

\move(104 10) \clvec (107 23)(149 23)(152 10)
\move(110 10) \clvec (112 13)(114 13)(116 10)
\move(122 10) \clvec (124 17)(138 17)(140 10)
\move(128 10) \clvec (131 21)(161 21)(164 10)
\move(146 10) \clvec (148 15)(156 15)(158 10)

\move(104 10) \fcir f:0 r:0.4
\move(110 10) \fcir f:0 r:0.4
\move(116 10) \fcir f:0 r:0.4
\move(122 10) \fcir f:0 r:0.4
\move(128 10) \fcir f:0 r:0.4
\move(134 10) \fcir f:0 r:0.4
\move(140 10) \fcir f:0 r:0.4
\move(146 10) \fcir f:0 r:0.4
\move(152 10) \fcir f:0 r:0.4
\move(158 10) \fcir f:0 r:0.4
\move(164 10) \fcir f:0 r:0.4

\move(105 7) \textref h:R v:C \htext{\scriptsize{$1^{(1)}$}}
\move(111.5 7) \textref h:R v:C \htext{\scriptsize{$1^{(2)}$}}
\move(118 7) \textref h:R v:C \htext{\scriptsize{$2^{(0)}$}}
\move(124.5 7) \textref h:R v:C \htext{\scriptsize{$2^{(1)}$}}
\move(131 7) \textref h:R v:C \htext{\scriptsize{$2^{(2)}$}}
\move(132 10) \textref h:R v:C \htext{\footnotesize{$*$}}
\move(136 7) \textref h:R v:C \htext{\scriptsize{$3^{(1)}$}}
\move(142 7) \textref h:R v:C \htext{\scriptsize{$4^{(0)}$}}
\move(144 10) \textref h:R v:C \htext{\footnotesize{$*$}}
\move(148.5 7) \textref h:R v:C \htext{\scriptsize{$5^{(1)}$}}
\move(155 7) \textref h:R v:C \htext{\scriptsize{$6^{(0)}$}}
\move(161.5 7) \textref h:R v:C \htext{\scriptsize{$7^{(0)}$}}
\move(168 7) \textref h:R v:C \htext{\scriptsize{$8^{(0)}$}}
}
\label{linear rep.-add marks}
\end{figure}

For a linked partition $\pi$ of $[n]$, the marked linear representation
of $\pi$ has $k$ marks if $\gpi$ has $k+1$ connected components  (as there is
no mark before the vertex $1$).
In every connected component of $\gpi$, the number of vertices is one greater
than the number of arcs, so in $\gpi$,
\begin{eqnarray*}
n &=&\#vertices\\
    &=&\#arcs+\#components\\
    &=&\#arcs+\#marks+1.
\end{eqnarray*}
For a linked cycle $\hat  \pi$ whose underlying  linked partition is $\pi$,
the marked linear representation of $\hat \pi$  has the same number of arcs and marks
as  that of  $\pi$. Hence
\begin{prop}\label{arc+mark}
Any marked linear representation of a linked cycle of [n] satisfies
$$\#arcs + \#marks = n-1.$$
\end{prop}

\textsc{Remark}.
It is clear that both the cycle representation $\ghpi^c$ and the linear
representation $\ghpi^l$ uniquely determine $\hat \pi$.
In the marked linear representation of $\hat \pi$, we may remove the
labeling on the vertices, as it can be recovered by the arcs and marks.
More precisely, suppose $G$ is a graph of a (partial) matching  whose vertices
are listed on a horizontal line, and where some vertices
$2, 3,\dots,n$  have marks before them.  If the total number of arcs
and marks is
$n-1$, then one can partition the vertices of $G$ into $n$ intervals
in the following way.
Start from the left-most vertex, end an interval before each righthand
endpoint or mark. Then for the $i$th interval, label the righthand
endpoint by $i^{(0)}$, the lefthand endpoints consecutively by  $i^{(1)},
i^{(2)}, \dots$, and isolated point by $i^{(1)}$.

\subsection{Enumeration of linked cycles}

Let $lc_n$ be the cardinality of  $LC(n)$, the set of
linked cycles on $[n]$.
It is easy to get  $lc_1=1, lc_2=2, lc_3=7, lc_4=37$. For any $\hat \pi \in
LC(n)$, let $s(\hp)$  be the number of singly covered minimal
elements in $\hp$.  $s(\hp)$ is also the number of connected components in
the cycle representation $\ghpi^c$.
We denote by $f(n,m)$ the number of linked
cycles in $LC(n)$ with $s( \hp)=m$.
Let $size(\hp)=\sum_{C_i} |C_i|$ where the sum is over
all cycles of $\hp$, and $double(\hp)$ be the number of doubly covered
elements of $\hp$.

\begin{prop}\label{thm_f_n1}
The numbers $f(n,m)$ satisfy the recurrence
$$
f(n,m)=(2(n-1)-m)f(n-1,m) + f(n-1, m-1)
$$
with initial  values $f(1,1)=1$, $f(n,0)=0$ and $f(n,m)=0$ if $n
<m$.
\end{prop}
\begin{proof} Clearly $f(1,1)=1$ and $f(n,m)=0$ if $n <m$. For any
linked cycle $\hp \in LC(n)$, $1$ is always a singly covered minimal
element. Hence $f(n,0)=0$.

Given $\hp \in LC(n)$, after removing the vertex $n$ and all edges incident
to $n$ in the cycle representation $\ghpi^c$, we get a linked
cycle on $[n-1]$.

Conversely,
 given $\hp'=\{B_1, B_2, ..., B_k\} \in LC(n-1)$, we can obtain a
linked cycle  $\hp$ of $[n]$ by joining the element $n$ in one of the
following mutually exclusive ways.
\begin{enumerate}
\item $\hp=\hp' \cup \{n\}$, that is, $\hp$ is obtained from
$\hp'$ by adding a singleton block $\{n\}$. In this case
$s(\hp)=s(\hp')+1$.
\item $\hp$ is obtained from $\hp'$ by inserting  $n$ into an
existing cycle $(a_1a_2...a_k)$ of $\hp'$. Since $n$ can be
inserted after any element $a_i$, there are $size(\hp')$ many such
formed $\hp$. For each of them, $s(\hp)=s(\hp')$.
\item $\hp$ is obtained from $\hp'$ by adding a cycle of the form
$(i,n)$. Such a constructed $\hp$ is a linked cycle if and only if  $i$
  is singly covered and not a minimal element of $\hp'$. There are
$n-1-s(\hp')-double(\hp')$ many choices for $i$. For each $\hp$
in this case, $s(\hp)=s(\hp')$.
\end{enumerate}

Combining these three cases, and noting that
$size(\hp')=n-1+double(\hp')$, we have
$size(\hp')+n-1-s(\hp')-double(\hp')=2(n-1)-s(\hp')$, which leads to
the desired recurrence
$$
f(n,m)=(2(n-1)-m)f(n-1,m) + f(n-1, m-1).
$$
\end{proof}

The initial values of $f(n,m)$ are
\begin{center}
\begin{tabular}{c|cccc}
$n \setminus k $ & 1 & 2  & 3 & 4  \\ \hline
1                & 1 &    &   &   \\
2                & 1 & 1  &   &   \\
3                & 3 & 3  & 1 &   \\
4                & 15& 15 & 6 & 1
\end{tabular}
\end{center}

The set of numbers $\{f(n,m): n \geq m >0 \}$  are the coefficients of
Bessel polynomials $y_n(x)$ (with exponents in decreasing order), and has been
studied by Riordan \cite{Riordan}, and by W. Lang as the signless
\emph{$k$-Stirling numbers of the second kind} with $k=-1$,
\cite{Lang2000}. From \cite[pp77]{Riordan}, we deduce that
\begin{eqnarray}\label{for_f(n,m)}
f(n,m)=\frac{(2n-m-1)!}{(m-1)!(n-m)!2^{n-m}} = {2n-m-1 \choose m-1} \frac{(2n-2m)!}{(n-m)!2^{n-m}}.
\end{eqnarray}
A combinatorial proof of \eqref{for_f(n,m)} is given at the end of this
subsection, using a bijection between linked cycles and certain set
partitions. Another definition for $f(n,m)$ is given by
the coefficients in the
expansion of the operator $(x^{-1}d_x)^n$, i.e.,
$$
(x^{-1}d_x)^n = \sum_{m=1}^n (-1)^{n-m} f(n,m) x^{m-2n} d_x^m,
\qquad n \in N.
$$
An extensive algebraic treatment based on this equation was given in \cite{Lang2000}.
The numbers $\{f(n,m)\}$ can be recorded in an infinite-dimensional lower
triangular matrix. In particular, $lc_n$, the cardinality of
$LC(n)$, is the $n$-th row sum of the matrix. The exponential
generating function of $lc_n$ is given in Formula (54) of
\cite{Lang2000} as
$$
\sum_{n=1}^\infty lc_n \frac{x^n}{n!} =\exp(1-\sqrt{1-2x}) -1.
$$
In the OEIS \cite{Sloane} $\{ lc_n\}$  is the sequence A001515,
where the
following  recurrence is given.
\begin{equation}
a_n=(2n-3)a_{n-1}+a_{n-2}
\end{equation}
Here we present a combinatorial proof based on the structure of linked cycles.

\begin{prop}\label{thm_f_n2}
The numbers $f(n,m)$ satisfy the recurrence
$$
f(n,m)=mf(n,m+1) + f(n-1, m-1)
$$
with initial  values $f(1,1)=1$, $f(n,0)=0$ and $f(n,m)=0$ if $n
<m$.
\end{prop}
\begin{proof}
Let
$$
A=\{\hp \in LC(n)~:~s(\hp)=m \text{ and } n \text{ is not a singleton block}\},
$$
and
$$
B=\{\hp \in LC(n)~:~s(\hp)=m+1\},
$$
Clearly $|A|=f(n,m)-f(n-1, m-1)$.
We will construct an $m$-to-one correspondence between the sets $A$ and
$B$.

Let $\hp \in LC(n)$ with cycle representation $\ghpi^c$.
First we describe an operation $\tau$, which decomposes a
non-singleton connected component of $\ghpi^c$ into two.
Let $G_1$ be a connected component with at least two vertices. Assume
$i_1$ is the minimal vertex of $G_1$, which lies in the root cycle
$C_1=(i_1i_2\dots i_t)$. Then $i_t$ must be at least
$2$. The operation $\tau$ removes the arcs $(i_1,i_2)$ and $(i_2,i_3)$,
and adds an arc $(i_1,i_3)$.

The inverse operation $\rho$ of $\tau$ merges two connected components $G_1$ and
$G_2$ of a linked cycle as
follows. Assume the minimal elements of $G_1$ and $G_2$ are $i_1$ and
$j_1$ respectively and $i_1 < j_1$. Then $i_1$ and $j_1$ must be singly covered.
Assume the root cycle of $G_1$ is $C_1=(i_1i_2\dots i_k)$. The
operation $\rho$ inserts $j_1$ into $C_1$ to get $(i_1j_1i_2\dots
i_k)$, and keeps all other cycles unchanged.

Now we can define the $m$-to-$1$ correspondence between the sets $A$ and
$B$. For any $\hp \in A$ the cycle representation $\ghpi^c$ has $m$
connected components, where $n$ is not an isolated point.
Applying the operation $\tau$ to the component containing the vertex $n$,
we get a linked cycle with $m+1$ connected components.
(See Figure \ref{LC_n->LC_{n-1}} for an illustration.)

Conversely, given any $\hp\in B$, there are $m+1$ connected
components  in $\ghpi^c$. We may merge the component
containing $n$ with any other component to get a linked cycle in
$A$.  There are $m$ choices for the other components, hence we
get an $m$-to-one correspondence.

\begin{figure}[h,t]
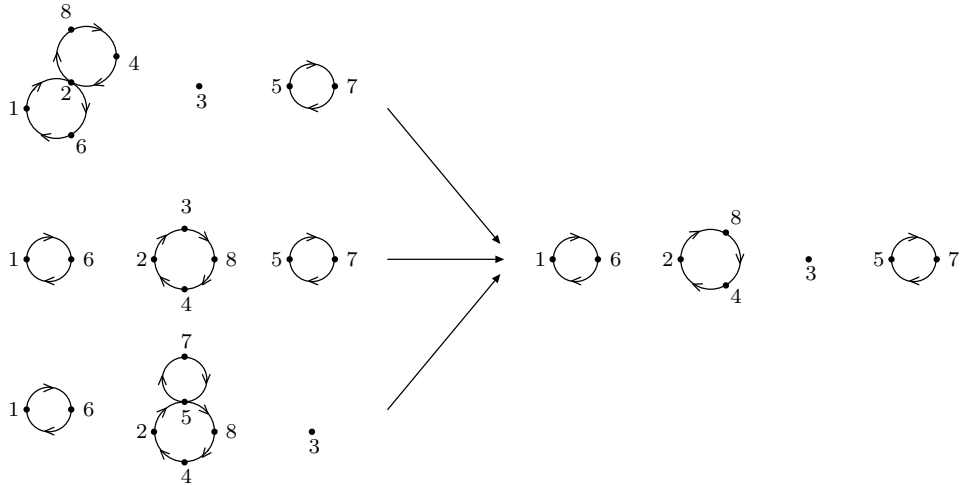

\centertexdraw{ \drawdim mm

\linewd 0.2

\move(6 50) \lcir r:4
\move(10 57) \lcir r:4

\move(2 50) \fcir f:0 r:0.4
\move(8 53.5) \fcir f:0 r:0.4
\move(8 46.5) \fcir f:0 r:0.4
\move(8 60.5) \fcir f:0 r:0.4
\move(14 57) \fcir f:0 r:0.4

\move(1 50) \textref h:R v:C \htext{\scriptsize{$1$}}
\move(8 52) \textref h:R v:C \htext{\scriptsize{$2$}}
\move(10 45) \textref h:R v:C \htext{\scriptsize{$6$}}
\move(8 63) \textref h:R v:C \htext{\scriptsize{$8$}}
\move(17 56) \textref h:R v:C \htext{\scriptsize{$4$}}

\arrowheadsize l:1.2 w:1
\move(3 52.5) \arrowheadtype t:v \avec(4 53.5)
\move(9.9 50.5) \arrowheadtype t:v \avec(9.9 49.5)
\move(5 46.2) \arrowheadtype t:v \avec(4 46.5)

\move(6.1 56.5) \arrowheadtype t:v \avec(6.1 57.5)
\move(11 60.8) \arrowheadtype t:v \avec(12 60.5)
\move(12 53.5) \arrowheadtype t:v \avec(11 53.2)

\move(25 53) \fcir f:0 r:0.4

\move(26 51) \textref h:R v:C \htext{\scriptsize{$3$}}

\move(40 53) \lcir r:3

\move(37 53) \fcir f:0 r:0.4
\move(43 53) \fcir f:0 r:0.4

\move(36 53) \textref h:R v:C \htext{\scriptsize{$5$}}
\move(46 53) \textref h:R v:C \htext{\scriptsize{$7$}}

\arrowheadsize l:1.2 w:1
\move(39.5 55.9) \arrowheadtype t:v \avec(40.5 55.9)
\move(40.5 50.1) \arrowheadtype t:v \avec(39.5 50.1)

\move(5 30) \lcir r:3

\move(2 30) \fcir f:0 r:0.4
\move(8 30) \fcir f:0 r:0.4

\move(1 30) \textref h:R v:C \htext{\scriptsize{$1$}}
\move(11 30) \textref h:R v:C \htext{\scriptsize{$6$}}

\arrowheadsize l:1.2 w:1
\move(4.5 32.9) \arrowheadtype t:v \avec(5.5 32.9)
\move(5.5 27.1) \arrowheadtype t:v \avec(4.5 27.1)

\move(23 30) \lcir r:4

\move(19 30) \fcir f:0 r:0.4
\move(27 30) \fcir f:0 r:0.4
\move(23 34) \fcir f:0 r:0.4
\move(23 26) \fcir f:0 r:0.4

\move(18 30) \textref h:R v:C \htext{\scriptsize{$2$}}
\move(30 30) \textref h:R v:C \htext{\scriptsize{$8$}}
\move(24 37) \textref h:R v:C \htext{\scriptsize{$3$}}
\move(24 24) \textref h:R v:C \htext{\scriptsize{$4$}}

\arrowheadsize l:1.2 w:1
\move(19.8 32.4) \arrowheadtype t:v \avec(20.6 33.2)
\move(25.4 33.2) \arrowheadtype t:v \avec(26.2 32.4)
\move(26.2 27.6) \arrowheadtype t:v \avec(25.4 26.6)
\move(20.6 26.8) \arrowheadtype t:v \avec(19.8 27.6)

\move(40 30) \lcir r:3

\move(37 30) \fcir f:0 r:0.4
\move(43 30) \fcir f:0 r:0.4

\move(36 30) \textref h:R v:C \htext{\scriptsize{$5$}}
\move(46 30) \textref h:R v:C \htext{\scriptsize{$7$}}

\arrowheadsize l:1.2 w:1
\move(39.5 32.9) \arrowheadtype t:v \avec(40.5 32.9)
\move(40.5 27.1) \arrowheadtype t:v \avec(39.5 27.1)

\move(5 10) \lcir r:3

\move(2 10) \fcir f:0 r:0.4
\move(8 10) \fcir f:0 r:0.4

\move(1 10) \textref h:R v:C \htext{\scriptsize{$1$}}
\move(11 10) \textref h:R v:C \htext{\scriptsize{$6$}}

\arrowheadsize l:1.2 w:1
\move(4.5 12.9) \arrowheadtype t:v \avec(5.5 12.9)
\move(5.5 7.1) \arrowheadtype t:v \avec(4.5 7.1)

\move(23 7) \lcir r:4
\move(23 14) \lcir r:3

\move(19 7) \fcir f:0 r:0.4
\move(27 7) \fcir f:0 r:0.4
\move(23 11) \fcir f:0 r:0.4
\move(23 17) \fcir f:0 r:0.4
\move(23 3) \fcir f:0 r:0.4

\move(18 7) \textref h:R v:C \htext{\scriptsize{$2$}}
\move(30 7) \textref h:R v:C \htext{\scriptsize{$8$}}
\move(24 9) \textref h:R v:C \htext{\scriptsize{$5$}}
\move(24 19) \textref h:R v:C \htext{\scriptsize{$7$}}
\move(24 1) \textref h:R v:C \htext{\scriptsize{$4$}}

\arrowheadsize l:1.2 w:1
\move(19.8 9.4) \arrowheadtype t:v \avec(20.6 10.2)
\move(25.4 10.2) \arrowheadtype t:v \avec(26.2 9.4)
\move(26.2 4.6) \arrowheadtype t:v \avec(25.4 3.6)
\move(20.6 3.8) \arrowheadtype t:v \avec(19.8 4.6)

\move(20.1 13.5) \arrowheadtype t:v \avec(20.1 14.5)
\move(25.9 14.5) \arrowheadtype t:v \avec(25.9 13.5)

\move(40 7) \fcir f:0 r:0.4

\move(41 5) \textref h:R v:C \htext{\scriptsize{$3$}}

\move(50 50) \arrowheadtype t:F \avec(65 32)
\move(50 30) \arrowheadtype t:F \avec(65.5 30)
\move(50 10) \arrowheadtype t:F \avec(65 28)

\move(75 30) \lcir r:3

\move(72 30) \fcir f:0 r:0.4
\move(78 30) \fcir f:0 r:0.4

\move(71 30) \textref h:R v:C \htext{\scriptsize{$1$}}
\move(81 30) \textref h:R v:C \htext{\scriptsize{$6$}}

\arrowheadsize l:1.2 w:1
\move(74.5 32.9) \arrowheadtype t:v \avec(75.5 32.9)
\move(75.5 27.1) \arrowheadtype t:v \avec(74.5 27.1)

\move(93 30) \lcir r:4

\move(89 30) \fcir f:0 r:0.4
\move(95 33.5) \fcir f:0 r:0.4
\move(95 26.5) \fcir f:0 r:0.4

\move(88 30) \textref h:R v:C \htext{\scriptsize{$2$}}
\move(97 35.5) \textref h:R v:C \htext{\scriptsize{$8$}}
\move(97 25) \textref h:R v:C \htext{\scriptsize{$4$}}

\arrowheadsize l:1.2 w:1
\move(90.6 33.3) \arrowheadtype t:v \avec(91.4 33.7)
\move(96.9 30.5) \arrowheadtype t:v \avec(96.9 29.5)
\move(91.4 26.3) \arrowheadtype t:v \avec(90.6 26.7)

\move(106 30) \fcir f:0 r:0.4

\move(107 28) \textref h:R v:C \htext{\scriptsize{$3$}}

\move(120 30) \lcir r:3

\move(117 30) \fcir f:0 r:0.4
\move(123 30) \fcir f:0 r:0.4

\move(116 30) \textref h:R v:C \htext{\scriptsize{$5$}}
\move(126 30) \textref h:R v:C \htext{\scriptsize{$7$}}

\arrowheadsize l:1.2 w:1
\move(119.5 32.9) \arrowheadtype t:v \avec(120.5 32.9)
\move(120.5 27.1) \arrowheadtype t:v \avec(119.5 27.1)

} \caption{The operation $\tau$ on three linked cycles with 3 components
 all gives the same linked cycle with 4 components.}\label{LC_n->LC_{n-1}}
\end{figure}
\end{proof}

\begin{thm}The numbers $lc_n$ $(n>1)$ satisfy the recursive relation
$$
lc_n=(2n-3)lc_{n-1}+lc_{n-2}
$$
with the initial values $lc_1=1$ and $lc_2=2$.
\end{thm}
\begin{proof}
\begin{eqnarray}
lc_n & = & \sum_{m=1}^n f(n,m) \nonumber \\
    & = & \sum_{m=1}^n [(2n-2-m)f(n-1,m) + f(n-1,m-1)]  \label{use_4.3}\\
    & = & \sum_{m=1}^n [(2n-2-m)f(n-1,m) +  (m-1)f(n-1,m) +
    f(n-2,m-2)] \label{use_4.4} \\
       & = & (2n-3)  \sum_{m=1}^n f(n-1, m) + \sum_{m=1}^n f(n-2, m-2)
    \nonumber \\
    & = & (2n-3)lc_{n-1} + lc_{n-2}, \nonumber
\end{eqnarray}
where the equation \eqref{use_4.3} is obtained by applying
Prop. \ref{thm_f_n1} to $f(n,m)$, and the
equation \eqref{use_4.4} is applying Prop. \ref{thm_f_n2} to
$f(n-1,m-1)$.
\end{proof}

In counting various kinds of set-partitions, Proctor
found a combinatorial interpretation for the
sequence $2,7,37,266,2431,\ldots$ as the number of
partitions of $[k]$ ($0\leq k\leq 2n$) into exactly $n$
blocks each having no more than  $2$ elements,
See \cite[\S7]{Proctor}. For example,  for $n=2$, there are
$7$ such partitions. They are
$p_1=\{\{1\},\{2\}\}$, $p_2=\{\{1,2\},\{3\}\}$,
$p_3=\{\{1,3\},\{2\}\}$, $p_4=\{\{1\},\{2,3\}\}$,
$p_5=\{\{1,2\},\{3,4\}\}$, $p_6=\{\{1,3\},\{2,4\}\}$,
$p_7=\{\{1,4\},\{2,3\}\}$. Denote by $P_2(n)$ the set of
such partitions for $0 \leq k \leq 2n$.

Our linked cycles provide  another combinatorial interpretation.
Using the marked linear representation,
we construct a bijection between the linked cycles on $[n+1]$
and the set  $P_2(n)$.

\noindent \underline{
{\bf A bijection $\gamma$ between  $LC(n+1)$ and $P_2(n)$.}}

Given a linked cycle $\hp \in LC(n+1)$ with the marked linear
representation,  removing the vertex labels and all isolated points,
then replacing each mark with a vertex,  and relabeling the vertices
by $1, 2, \dots, k$ from left to right, we get a graph of a partition
of  $[k]$ for some $0\leq k\leq 2n$, where each block has 1 or 2 elements.
By Prop. \ref{arc+mark}, there are exactly $n$ blocks in this partition.

Conversely, given a partition of $[k]$ ($0\leq k\leq 2n$) in $P_2(n)$
we represent it by a graph whose vertices are listed in a horizontal
line, and there is an arc connecting $i$ and $j$ if and only if
$\{i,j\}$ is a block. We can define a linked cycle on $[n+1]$ by
the following steps:
\begin{enumerate}
\item  Remove the labels of the vertices.
\item  Change each singleton block to a mark.
\item If a mark is followed by a righthand endpoint
or another mark,
add a vertex  right after it, and
\item If there is a mark before the first
vertex, add a vertex at the very beginning.
\end{enumerate}
The resulting graph is the marked linear representation of a linked
cycle of $[n+1]$. By the remark at the end of \S4.2, one can recover
the labeling of the vertices, and hence the linked cycle.

\begin{ex}
The linked cycle $\pi=(126)(284)(3)(57)$  corresponds to the
partition \\
$\{\{1,10\},\{2,3\},\{4,7\},\{5,12\},\{6\},\{8\},\{9,11\}\}$, where
$n=7$
and $7\leq k\leq 12$.\\
\\

\begin{figure}[h,t]
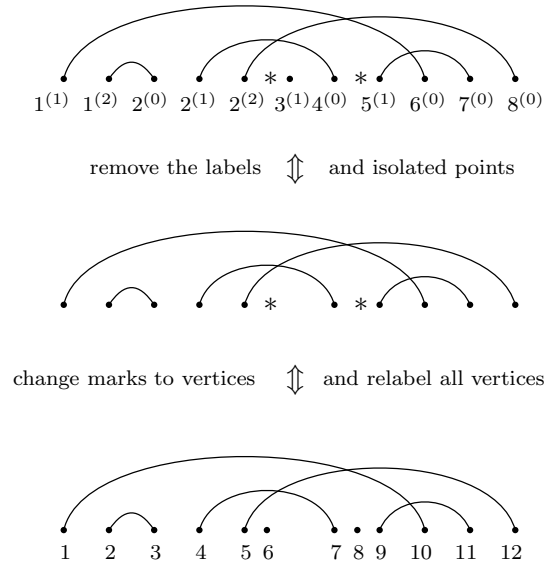

\centertexdraw{ \drawdim mm

\linewd 0.2

\move(4 70) \clvec (7 83)(49 83)(52 70)
\move(10 70) \clvec (12 73)(14 73)(16 70)
\move(22 70) \clvec (24 77)(38 77)(40 70)
\move(28 70) \clvec (31 81)(61 81)(64 70)
\move(46 70) \clvec (48 75)(56 75)(58 70)

\move(4 70) \fcir f:0 r:0.4
\move(10 70) \fcir f:0 r:0.4
\move(16 70) \fcir f:0 r:0.4
\move(22 70) \fcir f:0 r:0.4
\move(28 70) \fcir f:0 r:0.4
\move(34 70) \fcir f:0 r:0.4
\move(40 70) \fcir f:0 r:0.4
\move(46 70) \fcir f:0 r:0.4
\move(52 70) \fcir f:0 r:0.4
\move(58 70) \fcir f:0 r:0.4
\move(64 70) \fcir f:0 r:0.4

\move(5 67) \textref h:R v:C \htext{\scriptsize{$1^{(1)}$}}
\move(11.5 67) \textref h:R v:C \htext{\scriptsize{$1^{(2)}$}}
\move(18 67) \textref h:R v:C \htext{\scriptsize{$2^{(0)}$}}
\move(24.5 67) \textref h:R v:C \htext{\scriptsize{$2^{(1)}$}}
\move(31 67) \textref h:R v:C \htext{\scriptsize{$2^{(2)}$}}
\move(32.5 70) \textref h:R v:C \htext{$*$}
\move(37 67) \textref h:R v:C \htext{\scriptsize{$3^{(1)}$}}
\move(42 67) \textref h:R v:C \htext{\scriptsize{$4^{(0)}$}}
\move(44.5 70) \textref h:R v:C \htext{$*$}
\move(48.5 67) \textref h:R v:C \htext{\scriptsize{$5^{(1)}$}}
\move(55 67) \textref h:R v:C \htext{\scriptsize{$6^{(0)}$}}
\move(61.5 67) \textref h:R v:C \htext{\scriptsize{$7^{(0)}$}}
\move(68 67) \textref h:R v:C \htext{\scriptsize{$8^{(0)}$}}

\move(36 58) \textref h:R v:C \htext{$\Updownarrow$}

\move(64 58) \textref h:R v:C \htext{\scriptsize{remove the labels
    \qquad~ and isolated points}}

\move(4 40) \clvec (7 53)(49 53)(52 40)
\move(10 40) \clvec (12 43)(14 43)(16 40)
\move(22 40) \clvec (24 47)(38 47)(40 40)
\move(28 40) \clvec (31 51)(61 51)(64 40)
\move(46 40) \clvec (48 45)(56 45)(58 40)

\move(4 40) \fcir f:0 r:0.4
\move(10 40) \fcir f:0 r:0.4
\move(16 40) \fcir f:0 r:0.4
\move(22 40) \fcir f:0 r:0.4
\move(28 40) \fcir f:0 r:0.4
\move(40 40) \fcir f:0 r:0.4
\move(46 40) \fcir f:0 r:0.4
\move(52 40) \fcir f:0 r:0.4
\move(58 40) \fcir f:0 r:0.4
\move(64 40) \fcir f:0 r:0.4

\move(32.5 40) \textref h:R v:C \htext{$*$}
\move(44.5 40) \textref h:R v:C \htext{$*$}

\move(36 30) \textref h:R v:C \htext{$\Updownarrow$}

\move(68 30) \textref h:R v:C \htext{\scriptsize{change marks to
    vertices \qquad~ and relabel all vertices}}

\move(4 10) \clvec (7 23)(49 23)(52 10)
\move(10 10) \clvec (12 13)(14 13)(16 10)
\move(22 10) \clvec (24 17)(38 17)(40 10)
\move(28 10) \clvec (31 21)(61 21)(64 10)
\move(46 10) \clvec (48 15)(56 15)(58 10)

\move(4 10) \fcir f:0 r:0.4
\move(10 10) \fcir f:0 r:0.4
\move(16 10) \fcir f:0 r:0.4
\move(22 10) \fcir f:0 r:0.4
\move(28 10) \fcir f:0 r:0.4
\move(31 10) \fcir f:0 r:0.4
\move(40 10) \fcir f:0 r:0.4
\move(43 10) \fcir f:0 r:0.4
\move(46 10) \fcir f:0 r:0.4
\move(52 10) \fcir f:0 r:0.4
\move(58 10) \fcir f:0 r:0.4
\move(64 10) \fcir f:0 r:0.4

\move(5 7) \textref h:R v:C \htext{\scriptsize{$1$}}
\move(11 7) \textref h:R v:C \htext{\scriptsize{$2$}}
\move(17 7) \textref h:R v:C \htext{\scriptsize{$3$}}
\move(23 7) \textref h:R v:C \htext{\scriptsize{$4$}}
\move(29 7) \textref h:R v:C \htext{\scriptsize{$5$}}
\move(32 7) \textref h:R v:C \htext{\scriptsize{$6$}}
\move(41 7) \textref h:R v:C \htext{\scriptsize{$7$}}
\move(44 7) \textref h:R v:C \htext{\scriptsize{$8$}}
\move(47 7) \textref h:R v:C \htext{\scriptsize{$9$}}
\move(53 7) \textref h:R v:C \htext{\scriptsize{$10$}}
\move(59 7) \textref h:R v:C \htext{\scriptsize{$11$}}
\move(65 7) \textref h:R v:C \htext{\scriptsize{$12$}}

} \caption{The procedure from a linked cycle to a partition.
}\label{lc-bij.}
\end{figure}

\end{ex}

We conclude this subsection with a combinatorial proof of Formula
\eqref{for_f(n,m)}.  Recall that $f(n,m)$ is the number of linked
cycles  on $[n]$ with $m$  singly covered minimal elements.
Under the above bijection $\gamma$, it is the number of partitions
in $P_2(n-1)$ with $m-1$ isolated points, and hence $n-m$ blocks of
size $2$. For such a partition, the total number of points is $k=2(n-m)+m-1
=2n-m-1$. To obtain such a partition, we can first choose $m-1$
elements from $[k]=[2n-m-1]$ as isolated points, and then
construct a complete matching on the remaining $2(n-m)$ elements.
The number of ways to do this is then
$$
{2n-m-1 \choose m-1} \frac{2(n-m)!}{(n-m)! 2^{n-m}}.
$$

\subsection{Crossings and nestings of linked cycles}

In this subsection we present results on the enumeration of crossings and
nestings, as well as 2-crossings and 2-nestings, for linked cycles. This was
our  original motivation  to introduce the notion of linked cycles.

Given a linked cycle $\hp \in LC(n)$ with linear representation
$\ghpi^l$. Denote by $\ml=\ml(\pi)$ the vertex labeling of $\ghpi^l$. Two
different linked cycles $\hp$ and $\hp'$ on $[n]$ may have the same vertex
labeling. If this happens, then  $\hp$ and $\hp'$ share the following
properties:
\begin{enumerate}
\item $\hp$ and $\hp'$ have the same number of cycles.
\item $\hp$ and $\hp'$ have the same set of singly covered minimal
elements;
\item  $\hp$ and $\hp'$ have the same set of doubly covered elements;
\item Each cycle $C_i$ of $\hp$ can be paired with a unique cycle  $C_i'$  of
$\hp'$ such that $|C_i|=|C_i'|$, and $C_i$ and $C_i'$ have the same minimal
element.
\end{enumerate}
Fix a vertex labeling $\ml$, denote by
$LC_n(\ml)$ the set of all linked cycles $\hp$ with $\ml(\hp)=\ml$.
In $\ml$, if a vertex has a label $i^{(0)}$, we say it is a lefthand endpoint; if
it has
a label $i^{(k)}$ with $k \geq 1$, we say it is a righthand endpoint; if it has a
label $i$, we say it is an isolated point.
A $\hp \in LC_n(\ml)$ corresponds to a matching between the set of
lefthand endpoints to the set of righthand endpoints.

Let $k \geq 2$ be an integer. A $k$-crossing of $\hp$ is a set of $k$
arcs $(x_1,y_1)$, $(x_2, y_2), \dots, (x_k, y_k)$ of $\ghpi^l$ such that
the vertices appear in the order $x_1, x_2, \dots, x_k, y_1, y_2, \dots,
y_k$
from left to right. A $k$-nesting is a set of $k$ arcs $(x_1, y_1)$, $(x_2,
y_2), \dots,
(x_k, y_k)$ such that the vertices appear in the order $x_1, x_2, \dots,
x_k, y_k, \dots,
y_2, y_1$.     Denoted by $cr_k(\hp)$ the number of $k$-crossings of $\hp$,
and $ne_k(\hp)$ the number of $k$-nestings of $\hp$. Finally,
let $cr(\hp)$ be the maximal $i$ such that $\hp$ has a $i$-crossing
and $ne(\hp)$ the maximal $j$ such that $\hp$ has a $j$-nesting.

Our first result is an analog of Theorem \ref{thm_cr_ne},
on the joint generating function  of $cr_2$ and $ne_2$.
For any lefthand endpoint $i^{(0)}$, let
$$
h(i)=|\{\text{righthand endpoints on the right of $i^{(0)}$}\}|-
|\{\text{lefthand endpoints on the right of $i^{(0)}$}\}|.
$$
Then

\begin{thm} \label{2cr_lc}
\begin{eqnarray*}
& &  \sum_{\hp \in LC_n(\ml)} x^{cr_2(\hp)} y^{ne_2(\hp)} = \sum_{\hp \in
LC_n(\ml)} x^{ne_2(\hp)} y^{cr_2(\hp)} \\
& = & \prod_{i^{(0)} \in \ml} (x^{h(i)-1} + x^{h(i)-2}y+\dots
+x^{h(i)-k}y^{k-1} +\cdots +xy^{h(i)-2}+y^{h(i)}).
\end{eqnarray*}
In particular, the statistics $cr_2$ and $ne_2$ have a symmetric joint
distribution over each
set $LC_n(\ml)$.
\end{thm}
The proof is basically the same as that of Theorem \ref{thm_cr_ne}. It is even simpler
since every
vertex is the endpoint of at most one arc, that is, all $r_i=1$ in the proof
of Theorem \ref{thm_cr_ne}.

Perhaps more interesting is the joint distribution of $cr(\hp)$ and
$ne(\hp)$ over
$LC_n(\ml)$. Recall the following result in \cite{CDDSY}:
Given a partition $P$ of $[n]$, let
\begin{eqnarray*}
\min(P) & = & \text{\{minimal block elements of $P$\},}\\
\max(P) & = & \text{\{maximal block elements of $P$\}.}
\end{eqnarray*}
Fix $S,T\subset[n]$ with $|S|=|T|$. Let $P_n(S,T)$ be the set
$\{P\in\Pi_n: \min(P)=S, \max(P)=T,\}$. Then
\begin{thm}[CDDSY] \label{CDDSY-main}
$$
\sum_{P\in P_n(S,T)}x^{cr(P)}y^{ne(P)}
=\sum_{P\in P_n(S,T)}x^{ne(P)}y^{cr(P)}.
$$
That is, the statistics $cr(P)$ and $ne(P)$ have a symmetric joint
distribution over each set $P_n(S,T)$.
\end{thm}

Although the \emph{standard representation} for a partition of $[n]$
given in \cite{CDDSY}
is different than the linear representation defined in the present paper,
they coincide on partial matchings.
 View the graph $\ghpi^l$ as the graph of a partial matching $P$.
Observe that by fixing the vertex labeling,
we actually fix the number of vertices in $\ghpi^l$, the set of
minimal block elements of $P$ (which is the set of lefthand endpoints and
isolated points), and the set of maximal block element of $P$ (which is
the set of righthand endpoints and isolated points). Taking all linked cycles
with the vertex labeling $\ml$ is equivalent to taking all possible
partial matchings with the given sets of minimal block elements and
maximal block elements.
Hence Theorem \ref{CDDSY-main} applies to the set of linked cycles,
and we obtain the following theorem.

\begin{thm}
$$
\sum_{\hp \in LN_n(\ml)} x^{cr(\hp)} y^{ne(\hp)} = \sum_{\hp \in LN_n(\ml)}
x^{ne(\hp)} y^{cr(\hp)}.
$$
That is, the statistics $cr(P)$ and $ne(P)$ have a symmetric joint
distribution over each set $LC_n(\ml)$.
\end{thm}

\section*{Acknowledgments}

The authors thank Ken Dykema for  introducing to us  the notion
of noncrossing linked partitions.  We also thank Robert Proctor
for sharing a preprint of \cite{Proctor} with us.

\end{document}